\documentclass[]{article}
\usepackage{amsmath}
\usepackage[top=1in,bottom=1in,left=1in,right=1in,includefoot]{geometry}
\usepackage{graphicx}
\usepackage{color}
\usepackage{algorithmic} %for describing algorithms
\usepackage{pdflscape}
\usepackage{rotating}
\usepackage{epstopdf}
\usepackage{tabu}
\usepackage{chngcntr}
\usepackage{amsfonts}
\newcommand{\norm}[1]{\lVert#1\rVert}
\definecolor{cerulean}{rgb}{.05,0.5,.8}

\definecolor{darkgreen}{rgb}{0,0.5,0}
\definecolor{pink}{rgb}{.7,0,.7}

\definecolor{purple}{rgb}{.5,.1,.9}
\title{Classification of weighted networks through mesoscale homological features}
\author{Ann Sizemore, Chad Giusti, Danielle Bassett}
\date{}
\begin{document}
\maketitle

\begin{abstract}{As complex networks find applications in a growing range of disciplines, the diversity of naturally occurring and model networks being studied is exploding. The adoption of a well-developed collection of network taxonomies is a natural method for both organizing this data and understanding deeper relationships between networks. Most existing metrics for network structure rely on classical graph-theoretic measures, extracting characteristics primarily related to individual vertices or paths between them, and thus classify networks from the perspective of local features. Here, we describe an alternative approach to studying structure in networks that relies on an algebraic-topological metric called \emph{persistent homology}, which studies intrinsically mesoscale structures called \emph{cycles}, constructed from \emph{cliques} in the network. We present a classification of 14 commonly studied weighted network models into four groups or classes, and discuss the structural themes arising in each class. Finally, we compute the persistent homology of two real-world networks and one network constructed by a common dynamical systems model, and we compare the results with the three classes to obtain a better understanding of those networks.}

%%%% If classification number provided then

\end{abstract}

\renewcommand\thefigure{\arabic{figure}}
\section{Introduction}
\setcounter{figure}{-1}
\setcounter{table}{-1}

Driven by applications in fields as diverse as robotics, neuroscience, and economics, the quantity and complexity of available network models are growing rapidly. Certain structural themes, derived from classical graph-theoretic measures, are commonly used as organizing principles for understanding the general properties of and relationships between these models. Perhaps the best-known example of such a recurring theme is \emph{small-world} structure \cite{watts1998collective}, characterized by a combination of small characteristic path length and large clustering coefficient. Simply observing that a network model has this small-world property provides important information about the organization and behavior of the studied system. 

Yet, a single property like small-worldness clearly does not provide a complete summary of network structure, and many small-world networks exhibit very different behavior \cite{amaral2000classes}. Therefore, it is often useful to classify networks using multiple taxonomies, each providing a different lens into the characteristics of the system. For example, node connectivity patterns can be used to partition networks into classes with stringy-periphery and multi-star structure \cite{guimera2007classes}. Moreover, connectivity measures averaged over nodes can be used to separate certain real-world networks from others: for example, distinguishing biomolecular, citation, and computer networks from ecology, transportation, social, and communication networks \cite{kantarci2013classification}. Finally, mesoscale properties of a network, such as community structure, can be used to partition networks into \emph{similarity} or \emph{interaction} classes, reflecting the method of network construction \cite{onnela2012taxonomies}.  

\begin{figure}[h!]
\begin{centering}
 \includegraphics[width=6in]{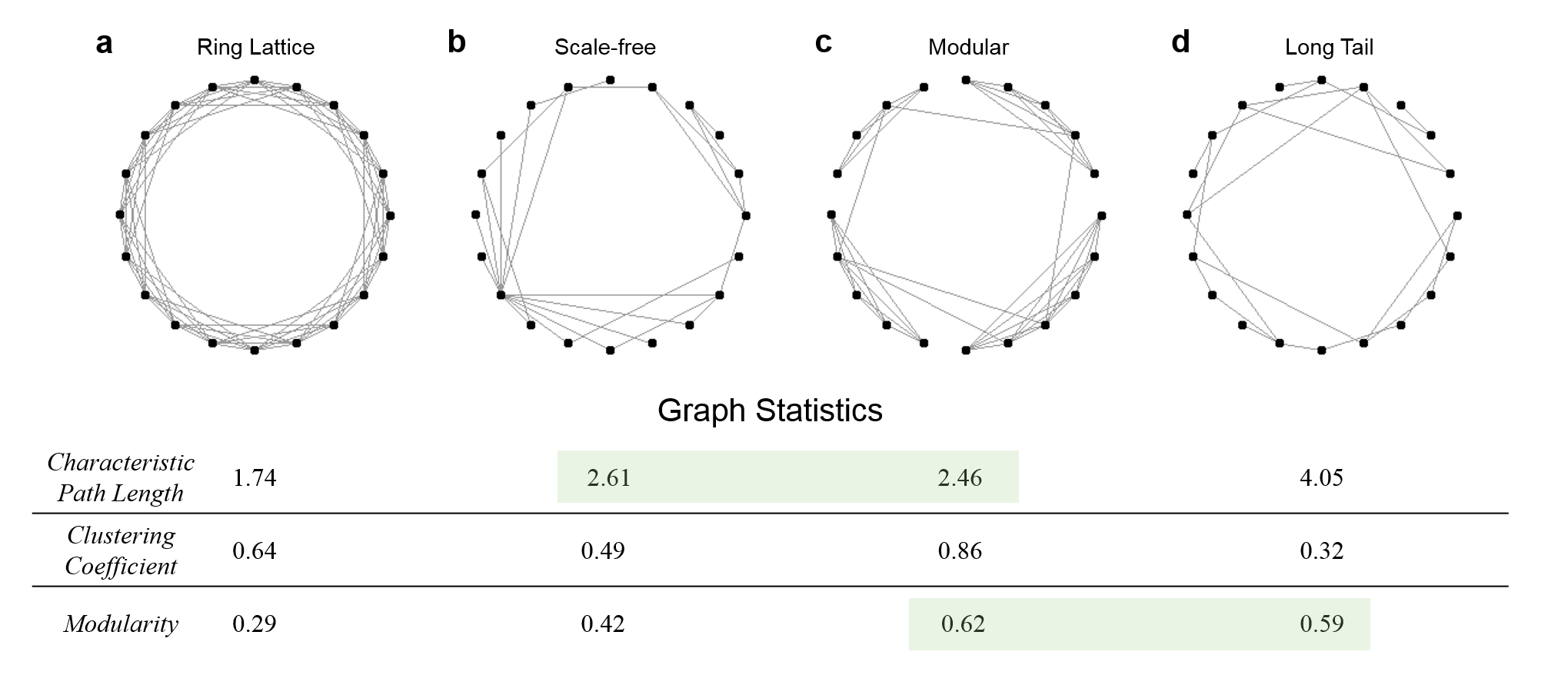}
\caption{Global Structure Provides an Alternative Perspective on Network Similarity. \emph{(a)} Regular ring lattice network \cite{watts1998collective}. \emph{(b)} Scale-free network with low average degree \cite{price1965networks}. \emph{(c)} Modular network of four components arranged cyclically \cite{gallai1967transitiv}. \emph{(d)} Connected circuit of ten nodes with a long tail. All networks are arranged to portray cyclic global structure. Green shaded regions contain very similar values within a particular graph statistic. Table below includes calculated graph statistics for the above network. Characteristic path length is the average of all shortest path lengths between nodes \cite{albert2002statistical}. The clustering coefficient measures the density of connections surrounding nodes \cite{luce1949method}. Modularity detects the separability of the network into clustered groups \cite{newman2006modularity}.}
\label{fig:rings} 
\end{centering}
\end{figure}

However, such classical graph-theoretic measures are usually \emph{local}, measuring properties of single vertices or of paths between fixed pairs, often then aggregating these measures to obtain a single network statistic. Such measures produce classifications that can be insensitive to higher order dependencies or structures in the network. For example, Fig.~\ref{fig:rings} shows four graphs on 20 nodes, which are considered quite different by classical graph statistics, yet carry a consistent global structure. The ring lattice in Fig.~\ref{fig:rings}a has identical statistics at each vertex; the scale-free network in Fig.~\ref{fig:rings}b has a low, but varying average degree; the network in Fig.~\ref{fig:rings}c has four modules with high internal degree, connected through hubs; and, the graph in Fig.~\ref{fig:rings}d has a high characteristic path length. The characteristic path length and modularity statistics find a pair of these four graphs to be very similar while the clustering coefficient values vary widely (Fig.~\ref{fig:rings} table). In spite of these local differences, however, all four networks have a consistent global structure in the form of a large-scale closed circuit.

Here, we describe a collection of statistics that detects such global structural features in a network, based on an algebraic-topological technique called \emph{homology} applied to the \emph{clique complex} of a graph. Recall that a \emph{clique} in a graph is an all-to-all connected subgraph; the \emph{clique complex} of a graph is a combinatorial object\footnote{An \emph{abstract simplicial complex}. For details, see \cite{ghrist2014elementary,ghrist2008barcodes,carlsson2009topology}.} whose constituent elements are the cliques of the graph. Homology measures how those cliques assemble to form particular loop patterns called \emph{cycles}, which correspond to our features of interest. For example, homology summarizes each network in Fig. \ref{fig:rings} as having one component and a single one-dimensional cycle (enclosing one two-dimensional void). Thus, through this lens these are four essentially identical graphs. 

Homology is, \emph{a priori}, a measurement of structure in unweighted graphs. However, a natural and powerful extension, \emph{persistent homology}, allows us to apply it as a measurement of weighted networks by tracking the evolution of cycles through the sequence of unweighted graphs obtained by thresholding the network at distinct edge weights (see Fig. \ref{fig:PH}, Methods, or \cite{carlsson2009topology} for more details). It is useful to think of such a sequence of graphs as obtained by beginning with a set of vertices and then adding new edges (and thus new cliques) as the threshold weight decreases. Thus, new cycles are \emph{born} and \emph{die} as they are filled in, providing a window into the global organization of the weights throughout the graph \cite{ghrist2008barcodes}. We use a variant of standard persistent homology that discards the threshold values, as doing so reduces sensitivity to choice of sampling distributions within models\footnote{For example, this makes the methods insensitive to potentially unknown geometric choices of scale.} \cite[Supplementary Information]{giusti2015clique}.

We illustrate the power of measuring mesoscale homological features using persistent homology in two ways. First, we show that this structure induces a classification of fourteen network models, including embedded and non-embedded, random and controlled networks, into four natural classes. The first class contains highly structured networks with a dense core, the second those with moderately sized clusters, the third random graphs subject to structural constraints, and the fourth class highly random networks. In addition, we compare real-world networks (created from brain imaging data and interactions of citric acid cycle enzymes) and dynamical systems networks (created from correlations between pairs of Kuramoto oscillators) to the identified classes to infer their mesoscale architecture. 

The paper is organized as follows: Methods includes a brief introduction to persistent homology, followed by a description of methods used in network generation and topological feature calculations; Results presents the resulting network classification and discusses the topological traits of networks in each class; Discussion interprets these results and suggests further applications.

%for comnet journal papers produced under \LaTeX\ using
%comnet.cls v1.5e.

\section{Methods}

\subsection{Homology and Applications}

We begin with a limited introduction to a persistent homology, describing only those concepts required to understand how our classification is performed. For a mathematical treatment and discussion of computational issues, we point the reader to the following useful references \cite{ghrist2014elementary,ghrist2008barcodes,zomorodian2005computing,carlsson2009topology}.

Let $G$ be an unweighted graph with $N$ vertices. A subset of $k$ vertices that are all-to-all connected is called a \emph{$k$-clique} of $G$. Recall that the convex hull of the $k$ points in general position is a $(k-1)$-dimensional region; the dimension, rather than the number of vertices, underlies the notion of \emph{degree} appearing in the literature, so we will use this indexing shift for consistency. We call any $m$-clique included in a $k$-clique ($m<k$) a \emph{face} and any clique that is not the face of another clique is termed \emph{maximal}. Denote the collection of all $k$-cliques in $G$ by $X_{k-1}(G)$, and let $X(G)=\{X_0(G),X_1(G),\dots,X_N(G)\}$ be the \emph{clique complex} of the network. Let $M_{k-1}(G)$ be the number of maximal $k$-cliques in $G$, and write the \emph{maximal clique vector} as $M(G) = \{M_0(G),M_1(G),\dots ,M_N(G)\}$ which records the number of maximal vertices ($1$-cliques), edges ($2$-cliques), triangles ($3$-cliques), etc. in a network. The clique number $\omega(G) = N+1$ records the number of vertices in the largest clique of $X(G)$. 

We define the \emph{boundary} of a given $k$-clique $\sigma$ as the set $\partial\sigma$ of all  $(k-1)$-faces of $\sigma$. The boundary of a set of $k$-cliques $\{\sigma_1,\sigma_2,\dots,\sigma_m\}$ is formed by taking the symmetric difference\footnote{The symmetric difference of sets $A$ and $B$, $A \Delta B = A\cup B \setminus A\cap B$. Recall that symmetric difference is an associative operation.} of the boundaries of the constituent cliques $$\partial\{\sigma_1,\sigma_2,\dots,\sigma_m\} = \partial\sigma_1 \;\Delta\; \partial\sigma_2\;\Delta\cdots\Delta\; \partial\sigma_m.$$ \noindent Geometrically, this corresponds to ``gluing together'' the cliques $\sigma_i$ along pairs of shared faces\footnote{The process we describe here is equivalent to computing homology with $\mathbb{Z}/2$ coefficients, the standard choice in topological data analysis.} to recover the ``shell'' encapsulating  $\{\sigma_1,\sigma_2, \dots \sigma_m\}$.

A \emph{($k-1$)-cycle} is a collection of $k$-cliques $\{\sigma_1,\sigma_2, \dots \sigma_m\}$ with $\partial\{\sigma_1,\sigma_2,\dots,\sigma_m\}=\emptyset$. Observe that any collection of $k$-cliques that appears as the boundary of a collection of $(k+1)$-cliques must be a $(k-1)$-cycle. However, a $(k-1)$-cycle may not be a boundary of higher dimensional cliques, and these cycles encode interesting, non-local structural relations about the arrangement of cliques in the underlying graph. 

 Two $k$-cycles $\ell_1 =\{ \sigma_{1,1},\sigma_{1,2},\sigma_{1,3},\dots , \sigma_{1,m}\}$ and $\ell_2=\{ \sigma_{2,1},\sigma_{2,2},\sigma_{2,3},\dots ,\sigma_{2,n}\}$  are \emph{equivalent} if $\ell_1 \Delta \ell_2$ is a boundary of a collection of higher dimensional cliques\footnote{Thus, in particular, any cycle that appears as a boundary is equivalent to an ``empty'' trivial cycle.}. The \emph{homology} of a clique complex is an enumeration of these equivalence classes of its cycles. By abuse, it is customary to refer to equivalence classes of $k$-cycles as $k$-cycles, and we will adopt that convention here. Denote the number of $k$-cycles in the clique complex of $G$ by $\beta_{k}(G)$ and record these in the Betti Sequence $\beta(G) = \{\beta_0(G), \beta_1(G),\dots\,\beta_N(G)\}$. This vector summarizes the non-bounding cycles found in the clique complex of the graph\footnote{Observe that $\beta_0(G)$ is the number of connected components in $G$: every vertex is a $0$-cycle, and any two such are equivalent if they share an edge.}.  Together, $M(G)$ and $\beta(G)$ provide a picture of the mesoscale structure of the unweighted graph $G$.

\subsection{Extension to Weighted Networks}

\begin{figure}[h!]
\begin{centering}
 \includegraphics[width=6in]{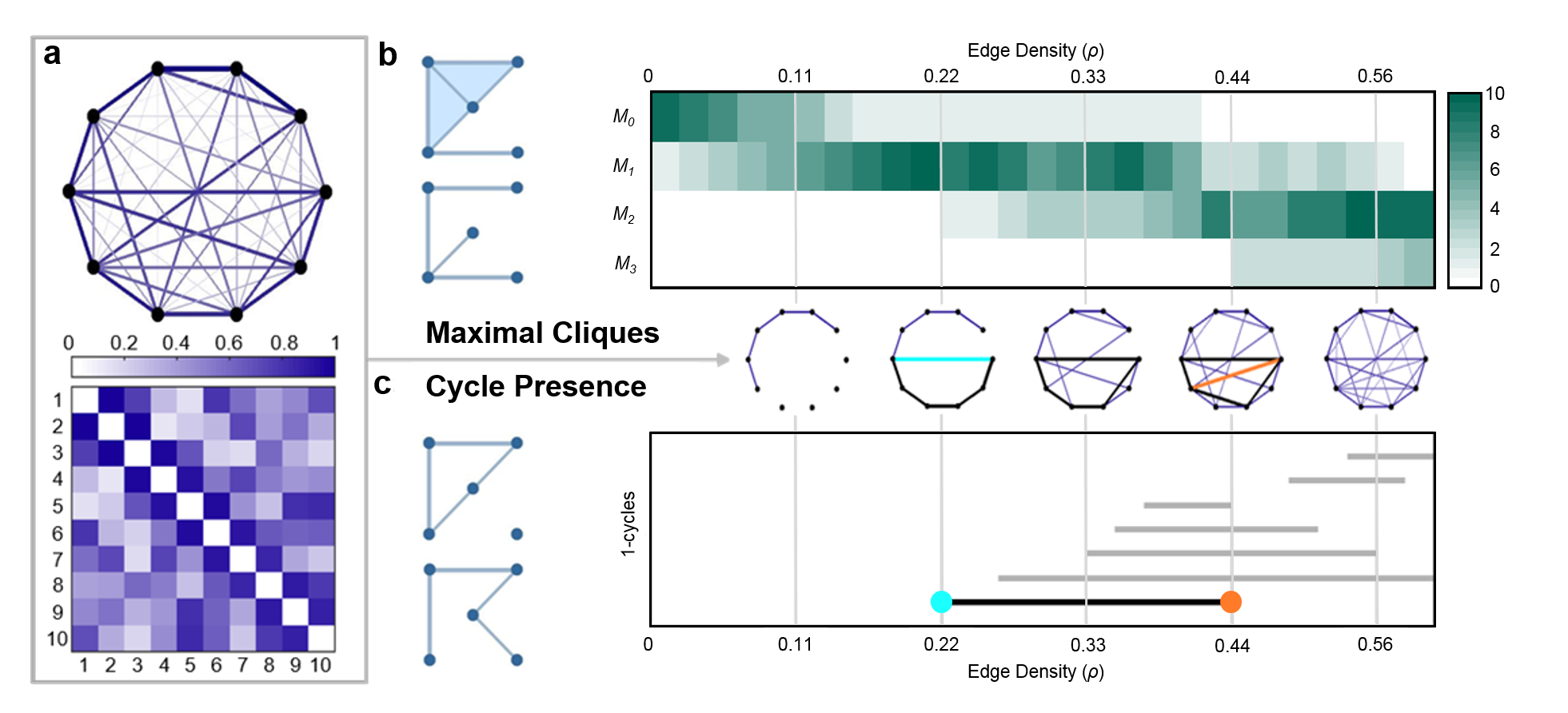}
\caption{Persistent Homology Detects Structure through Cliques and Cycles. \emph{(a)} (Top) Pictorial depiction of a weighted network on 10 nodes. Edge weights are indicated by darker color and line thickness. (Bottom) Adjacency matrix representation of the weighted graph above. \emph{(b)} (Left) Simple acyclic graphs with differing maximal clique distribution. (Right) Heat map showing the growth of $M_0$, $M_1$, $M_2$, and $M_3$ across edge density ($\rho$). Graphs below the horizontal axis of this heat map depict the network in Panel \emph{(a)} thresholded at the indicated edge density. \emph{(c)} (Left) Example graphs -- one cyclic and one acyclic -- but with identical numbers of cliques. (Right) Barcode plot showing the number of 1-cycles as a function of edge density. The black bar corresponds to the cycle created by the black edges in the filtration above the plot. This cycle is born with the addition of the cyan edge, shortened by subsequent edges, and finally killed by the orange edge. The birth and death density are indicated by the cyan and orange dots on the barcode plot, respectively.}
\label{fig:PH} 
\end{centering}
\end{figure}

Real world interactions are rarely binary, thus we require a translation of these measures to weighted networks in order to study empirical data. A weighted network (Fig. \ref{fig:PH}a, top) can be represented through its real-valued adjacency matrix (Fig. \ref{fig:PH}a, bottom). %, with entry $A_{i,j}$ the weight of the edge between nodes $i$ and $j$ 
The ordering on the edge weights induces a natural ordering on the edges, from strongest to weakest. We construct a sequence of unweighted graphs called a \emph{filtration} by beginning with an empty graph, and adding a single edge at a time per this ordering (Fig. \ref{fig:PH}, right); clearly each graph in the filtration is a subgraph of the next in a canonical fashion. 

For each graph in this progression, we compute the homology and record the maximal clique distribution and Betti sequence $\beta(G)$, all defined earlier in this Methods section.  Tracking the evolution of maximal cliques as the threshold drops provides a picture of the tightly-bonded structures in the network. However, while cliques only grow in size as the number of edges increase, homology behaves quite differently: cycles instead are ``born'', change in form, and eventually become boundaries of larger cliques and are ``killed''. 

As an example, the cycle highlighted in black (Fig. \ref{fig:PH}, right, sequence of graphs) formed from six edges is born when the cyan edge is added. Found on the \emph{barcode plot} below (Fig. \ref{fig:PH}c, right), the cyan dot indicates this edge density, termed the \emph{birth density}. The addition of edges at subsequent densities shortens the black cycle, until finally the addition of the orange edge divides the interior of the cycle (now four edges in length) into two 3-cliques. This kills the cycle, as it is now a boundary of 3-cliques. The orange dot on the barcode plot marks this edge density, the \emph{death density}. Cycle lifetime, indicated by the black bar (Fig. \ref{fig:PH}c, right), is the birth density subtracted from the death density, and will be used later to compute topological statistics (Section 2.3). Heuristically, cycles with longer lifetimes are called \emph{persistent}, as they must survive many edge additions, and are often considered the most essential topological features of the network. The homology of the clique complexes of all graphs in the filtration, along with this birth and death data, is referred to as the \emph{persistent homology} of the filtration.

%%%%%%%%%%%%%%%%%%%%%%%%%%%%%%%%%%%%%% NETWORK CONSTRUCTION%%%%%%%%%%%%%%%%%%%%%%%%%%%%%

\subsection{Network Construction}

\begin{figure}[!h]
\begin{centering}
 \includegraphics[width=6in]{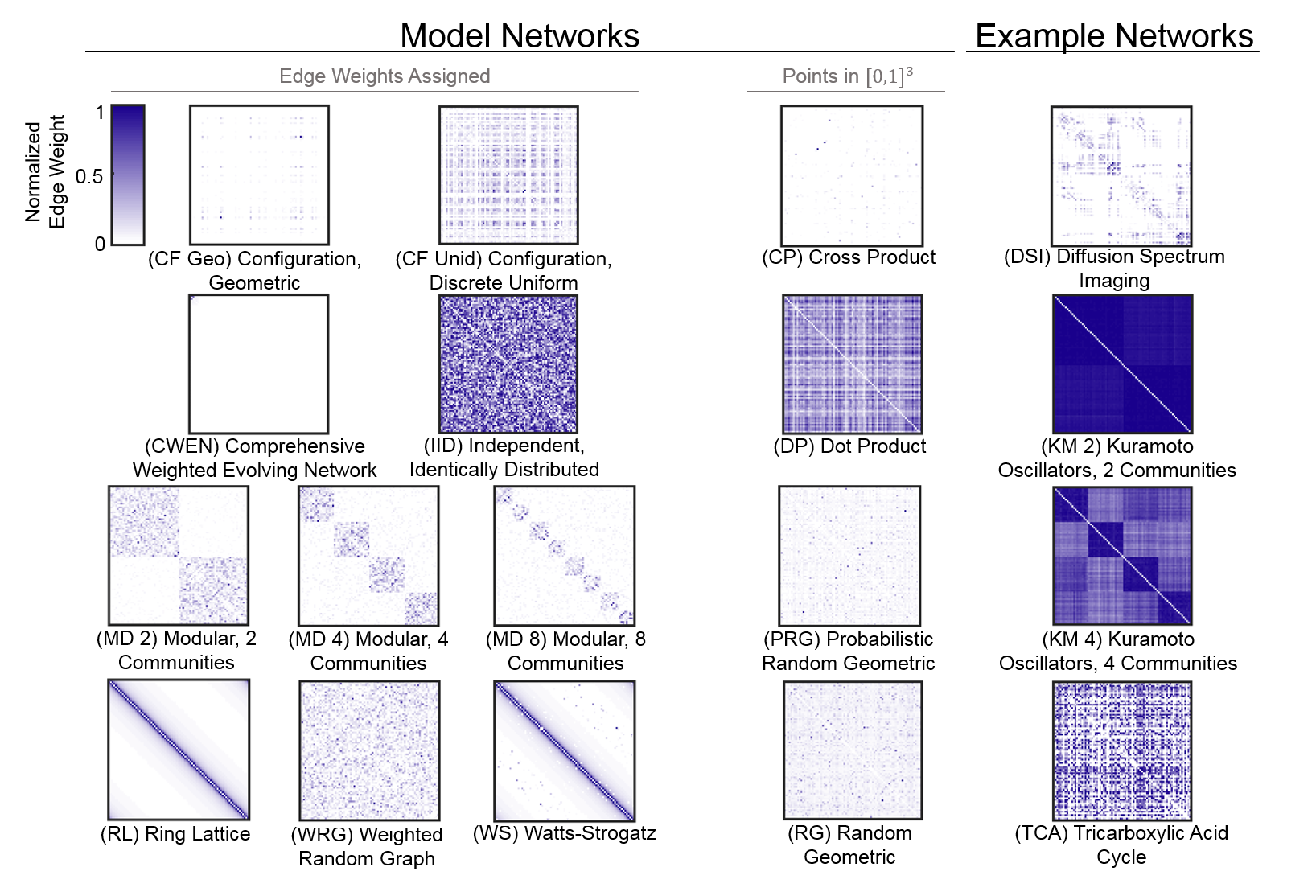}
\caption{Model and Example Networks. An example adjacency matrix of each model network is shown, identified with the network abbreviation found in Section 2.2. Edge weights in each network were normalized, and the resulting normalized weight is indicated by color.}
\label{fig:Mats} 
\end{centering}
\end{figure}

For this analysis, we  constructed model networks on 83 nodes\footnote{The features we study are essentially stable when normalized to scale for networks of this size and larger \cite[Supplementary Information]{giusti2015clique}.} (Fig. \ref{fig:Mats}) chosen for consistency with one of the real-world networks we later examine: a network of bundles of neuronal axons connecting large-scale regions of the human brain. We created models either by strictly identifying edge weights or through calculations on points in 3-space. We sampled each network model 100 times, and to ensure edge uniqueness in networks with redundant edge weights\footnote{Generically, edge weights are unique for ease of computation and comparison with real data. However, our computational methods do work for arbitrary networks.}, independent random noise was added from a uniform distribution on the interval $[0,0.001]$. % with $\lambda= 0.0005$, $\sigma =8.3\times10^{-8}$ after initial network construction, chosen to be smaller than the smallest difference between distinct weights in any model. This nose does not substantially impact our chosen measures as persistent homology is robust to a small amount of noise \cite{carlsson2009topology}. 
% cite cohen2007stability for stability of barcodes

We give a brief account of the network models here; detailed descriptions can be found in the Appendix. Code for generating all network models and computing all network statistics can be found in \cite{NGAT}.  
	\subsubsection{Networks from Edge Weighting Schemes}  

We study model networks that have been previously defined in the literature, particularly in comparisons to real-world networks in biology \cite{klimm2014resolving,lohse2014resolving}. The first group of these model networks that we study can be constructed from algorithms that define only edge weights (in the next section, we study networks that can be constructed from algorithms that assign locations of points in 3-space). Many of these algorithms produce a complete graph. For algorithms that could not produce a complete graph without sacrificing network traits, we chose parameter values to achieve an edge density of $\rho \sim 0.75$ to ensure capturing of third dimensional homology of the network, generically seen before $\rho = 0.6$ in random network models \cite{kahle2014topology}. Parameter values and further details can be found in \cite{NGAT}.\\
	\noindent Specifically, we tested the previously reported models:
	\begin{itemize}
	\item \textit{(CF) configuration model} with node strengths chosen from a discrete uniform (\textit{Unid}) distribution in the interval $[0,1000]$ or a  geometric (\textit{Geo}) distribution with $p=0.1$ \cite{serrano2005weighted}.
	\item \textit{(CWEN) comprehensive weighted evolving network} \cite{li2004comprehensive} with parameters chosen to create a scale-free strength distribution ($Pr(k) \sim k^{-\gamma}$) with $\gamma \sim 3$.
	\item \textit{(IID) independent and identically distributed} created by assigning edge weights random numbers from a uniform distribution on [0,1] \cite{kahle2011random}.
	\item \textit{(MD) modular} with varying numbers of communities. We created a binary modular graph using the Brain Connectivity Toolbox \cite{rubinov2010complex}, then we chose edge weights from a geometric distribution with probabilities based on whether endpoints of the edge lie in a single module or in distinct modules.
	\item \textit{(RL) ring lattice} \cite{watts1998collective} with edge weights inversely proportionate to hop distance along the perimeter of the lattice.
	\item \textit{(WRG) weighted random graph} created by determining edge weights from a geometric distribution \cite{garlaschelli2009weighted}.
	\item \textit{(WS) Watts-Strogatz} with edge switching parameter chosen to \cite{watts1998collective} maximize small-world propensity \cite{muldoon2015small}. 	
	\end{itemize}

	\subsubsection{Networks Created from Points in 3-Space} 

In addition to the model networks described above, we also examined model networks that were created from points in 3-space \cite{klimm2014resolving,lohse2014resolving}. Our motivation for studying this class of \emph{spatially embedded} networks \cite{barthelemy2011spatial} lies in the fact that in many real world systems, geometric constraints can govern network structure. We therefore tested the following network models, which each began with choosing 83 points uniformly at random in the $[0,1]^3$. Complete weighted networks were constructed from these points by performing the indicated calculations, additionally creating a correspondence of nodes and chosen points. 
	\begin{itemize}
	\item \textit{(CP) cross-product} For all $i, j$, $A_{i,j} = 1/\norm{\vec{i}\times \vec{j}}$ .
	
	\item \textit{(DP) dot product} For all $i,j$, $A_{i,j} = \vec{i} \cdot \vec{j}$. At each thresholded level, this is similar to the binary dot product graph designed by \cite{nickel2007random}.
	
	\item \textit{(PRG) probabilistic random geometric}: % A note to Chad and myself - this is actually called PRG_WS in our calculations, as my original PRG was pretty lame.
	With probability $p$, we swapped edges of a random geometric network similar to that of a Watts-Strogatz rewiring process on a ring lattice \cite{watts1998collective}. We chose $p$ to maximize small-world propensity \cite{muldoon2015small}. 
	
	\item \textit{(RG) random geometric}: The 83 nodes were wired minimally, assigning edge weights to $\frac{1}{d(\vec{i},\vec{j})}$ where $d(\vec{i},\vec{j})$ is the Euclidean distance between points $\vec{i}$ and $\vec{j}$ \cite{kahle2011random}. 
	
	\end{itemize}
	\subsubsection{Example Networks} 

Finally, below we consider two real-world networks and one network constructed by a common dynamical systems model.
	
	\begin{itemize}
\item \textit{(TCA) protein-protein interaction network of Tricarboxylic Acid Cycle Enzymes}: Also known as the Citric Acid cycle or Krebs cycle, the tricarboxylic acid cycle is a circular pathway that recovers energy from products of multiple metabolism channels. In \textit{Homo sapiens}, the 8 steps of the cycle are catalyzed by products of 17 genes \cite{TCA}. We used data from the STRING database to construct the protein-protein interaction network of these 17 gene products and their 66 strongest interacting partners (thus achieving an 83 node network). The STRING database includes the combined score, an overall rating of the certainty of interactions between any two proteins, which became the edge weights of the network \cite{szklarczyk2014string}.
	
	\item \textit{(DSI) structural human brain network constructed from Diffusion Spectrum Imaging}: Diffusion spectrum imaging (DSI) reveals the direction of water diffusion, from which a map of axonal trajectories between brain regions can be inferred. DSI data from an average of eight healthy adult humans \cite{gu2015controllability} induced a weighted network on 83 nodes according to the Lausanne scale33 parcellation of the brain into 83 regions of interest \cite{cammoun2012mapping}. We assigned the number of white matter tracts between regions $i,j$ to the edge weight between nodes $i,j$ \cite{klimm2014resolving,cieslak2014local}.
	
	\item \textit{(KM) Kuramoto Oscillator network}: This model of coupled oscillators has been used heavily in neuroscience to model neuronal behavior \cite{cumin2007generalising,breakspear2010generative}. As in \cite{bassett2014crosslink}, we built a network of 83 oscillators, assigning edge weights to the average correlation between oscillator pairs for networks with two (\textit{KM 2}) and four (\textit{KM 4}) communities. 
\end{itemize}

\subsection{Computations of Graph Statistics and Topological Features}
\textit{Graph Statistics}
We computed five representative graph statistics common in network analysis: the clustering coefficient ($C$), characteristic path length ($L$), local efficiency ($E_{loc}$), global efficiency ($E_{glob}$), and modularity ($Q$) of each network using the Brain Connectivity Toolbox \cite{rubinov2010complex}. See the Appendix for exact formulas.

\noindent\textit{Persistent Homology Computations}
Persistent homology allows us to computationally follow cliques and cycles as the edge density ($\rho$) increases. We calculated persistence intervals out to $\rho = 1$ with methods from \cite{henselmannovel} and used functions from the clique-top-master package \cite{giusti2015clique} for additional topological computations. 

\noindent\textit{Topological Statistics ($\overline{\beta}_{d}$,$\overline{\mu}_d$)}
We also computed statistics on the features recovered from persistent homology. In a given dimension $d$, we define $\overline{\beta}_{d}$, as the sum of the cycle lifetimes: let $x_m$ and $y_m$ be the birth and death densities of cycle $m$ in dimension $d$, respectively,

\begin{align}\label{bettibar1}
\begin{split}
\overline{\beta}_d = \sum_{m} (y_m - x_m)
\end{split}
\end{align}
summing over all $d$-cycles.

We also considered $\overline{\mu}_{d}$, the sum of the lifetimes weighted by the birth densities for each dimension $d$,
\begin{align}\label{bettibar2}
\begin{split}
\overline{\mu}_d = \sum_{m}x_m (y_m - x_m)
\end{split}
\end{align}
which is more sensitive to cycles with larger birth densities following \cite{adcock2012ring}. For dimension $0$, all cycles begin at edge density $0$. Thus, we assign $1/{83\choose 2}$ to the birth density of $0$-cycles for this computation.  Note all lifetimes fall within $[0,1]$.

\noindent\textit{Maximal Clique Distribution Tracking}
At each threshold, we extracted the number of maximal $k$-cliques in each dimension and estimated parameters $\mu$, $\sigma$, of the logarithmic normal distribution $f(x|\mu , \sigma) = \frac{1}{x\sigma \sqrt{2\pi}}\exp\{\frac{-(\ln x-\mu)^2}{2\sigma^2}\}$ by letting $\mu = \text{mean}(\ln(M_k))$ and $\sigma = \text{std}(\ln(M_k))$ \cite{mood1974introduction}. Due to computational limits, we recorded the maximal clique distribution out to $\rho =0.25$ for all networks. 

\noindent\textit{Betti curves ($\beta_d$)}
We recorded the number of cycles at each edge density, allowing us to see the fluctuation in cycle number as $\rho$ increases. This cycle number sequence in a particular dimension $d$ is the Betti curve $\beta_d$, and $\overline{\beta}_d = \int{\beta_d d\rho}$.

\subsection{Hierarchical Clustering}
To determine the structural similarities between networks based on their homological features, we perform an agglomerative hierarchical clustering on the models. In this method, each network begins as its own group, and the distance allowed between network features within one cluster is zero. As we allow larger distances between network features within a cluster, groups begin to merge until the number of desired clusters is reached \cite{maimon2005data}. We used Euclidean distances and chose $\overline{\beta}_d$,$\overline{\mu}_d$ values in dimensions 0-3 and the parameters $\mu$ and $\sigma$ of the maximal clique distribution logarithmic normal approximation for $k = 1,\dots ,\omega(G)$ as features. To determine the proper number of clusters, we analyzed silhouette plots shown in Fig. \ref{fig:silplot} \cite{kaufman1990finding}.

%%%%%%%%%%%%%%%%%% RESULTS %%%%%%%%%%%%%%%%%%%%%%%%%%%%%%%%%%%%%%

\section{Results}

	Using homological features, we are interested in determining networks with similar mesoscale structure. Persistent homology records such features by tracking cliques and cycles, which we use to cluster 14 model networks. 
		
\begin{figure}[h!]
			\centering
				\includegraphics[width= 6in]{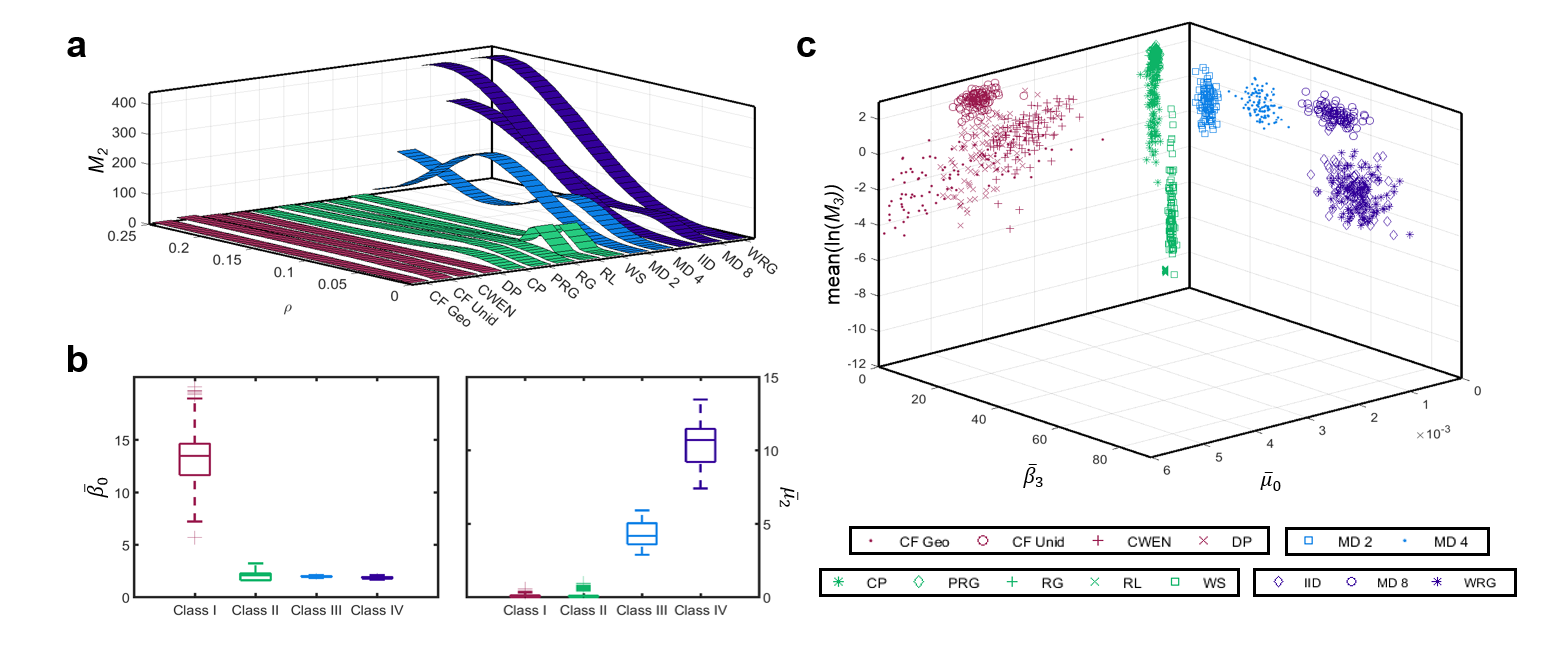}
				\caption{Topological Features Distinguish Novel Network Classes. Throughout the figure, the color of the symbol or line represents the results of the clustering. Class I shown in red, Class II in green, Class III in blue, and Class IV in purple. \emph{(a)} The value of $M_2$ tracked along edge density ($\rho$). Ribbon color is based on the clustering results. \emph{(b)} (Left) $\overline{\beta}_0$ and (Right) $\overline{\mu}_2$ results by class, indicated by box outline color. \emph{(c)} Identified network clusters plotted using three topological features $\overline{\mu}_0$, $\overline{\beta}_3$, and the average ln($M_3$) for $\rho \leq 0.25$. }
			\label{fig:Cluster}
		\end{figure}

		Clustering using Betti bar values and parameters from the maximal clique distribution fit partitions the model networks into four classes, shown in Fig. \ref{fig:Cluster} (see Fig. \ref{fig:silplot} for silhouette plots).  Recovered classes I-IV are distinguished by color: Class I shown in red, Class II in green, Class III in blue, and Class IV in purple. As an example of differences in maximal clique distribution across edge density ($\rho$), Fig. \ref{fig:Cluster}a shows $M_2$ (the number of maximal 3-cliques) across edge densities. Class IV members exhibit the highest $M_2$, followed by networks in Class III. Class II networks instead experience a short period with a high $M_2$, while Class I networks show a near zero $M_2$ for $\rho \leq 0.25$. Additionally, Fig. \ref{fig:Cluster}b shows the difference in $\overline{\beta}_{0}$ and $\overline{\mu}_2$ across recovered classes. Class I displays the highest $\overline{\beta}_0$ values, while Classes III and IV have distinctly higher $\overline{\mu}_2$ values. The scatter plot in Fig. \ref{fig:Cluster}c further illustrates class differences in homological features $\overline{\mu}_0$, $\overline{\beta}_3$, and the average ln($M_3$) for $\rho \leq 0.25$. Refer to Fig. \ref{fig:Sbbar1},\ref{fig:Sbbar2} for all $\overline{\beta}_d$ and $\overline{\mu}_d$ results as well as Fig. \ref{fig:sfit} for example logarithmic normal approximations of the maximal clique vector. These results demonstrate that the 14 model networks can be separated into distinct classes based on their topological properties.
		
\begin{figure}[h!]
			\centering
				\includegraphics[width=6in]{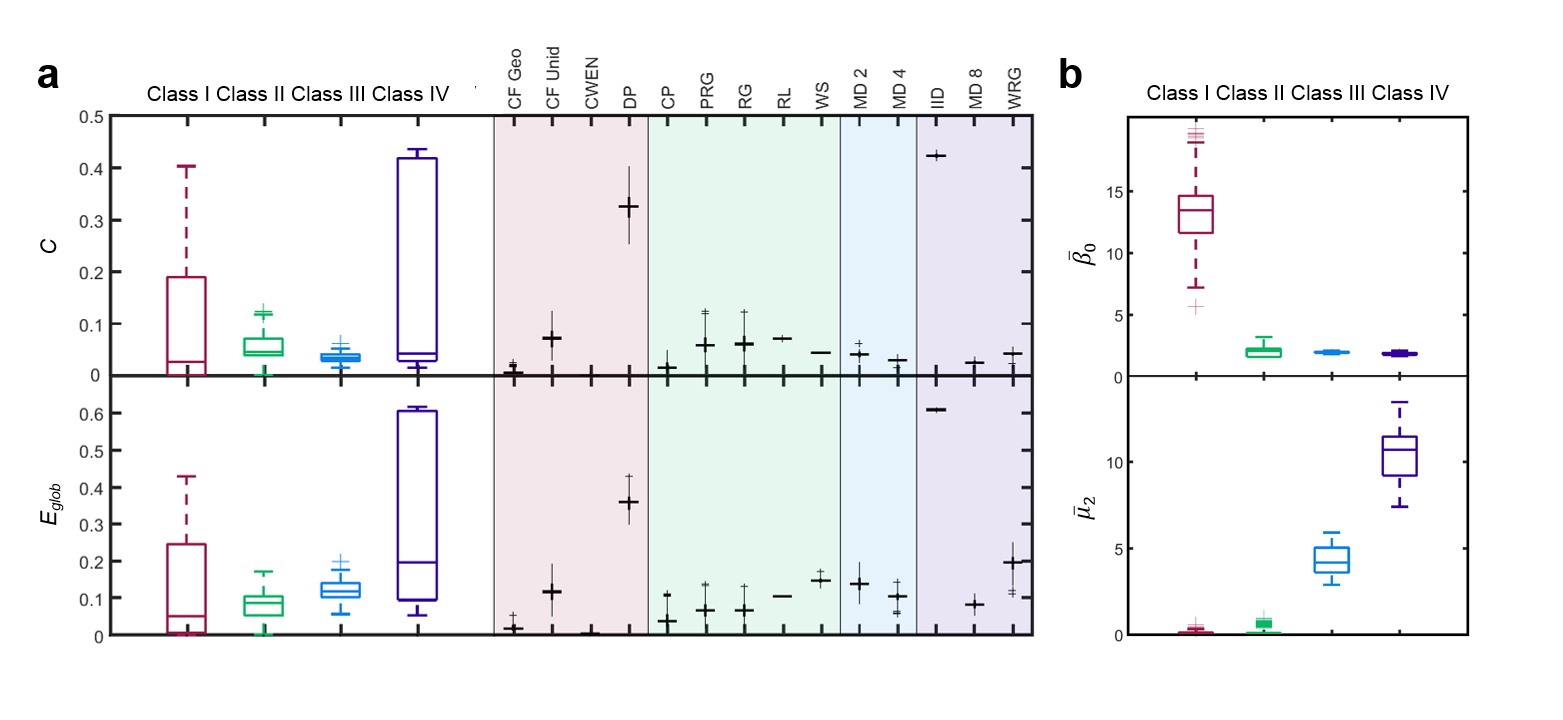}
				\caption{Common Graph Statistics are Insensitive to Topologically Defined Classes. \emph{(a)} (Top) Clustering coefficient ($C$) and (bottom) global efficiency ($E_{glob}$) shown by group and individually. Individual network boxes shown in black, shaded by group color. \emph{(b)} (Top) Calculated $\overline{\beta}_{0}$ and (bottom) $\overline{\mu}_{2}$ values shown by group (repeated from Fig. \ref{fig:Cluster}b).}
				\label{fig:GS} 
			
		\end{figure}

		Given the simple classification of model networks into four classes, provided by the topological characteristics, it is interesting to ask whether such insights could have been obtained from standard graph statistics. To address this question, we computed five graph statistics for model networks: clustering coefficient ($C$), global efficiency ($E_{glob}$), local efficiency ($E_{loc}$), characteristic path length ($L$), and modularity ($Q$). Box plots in Fig. \ref{fig:GS} show the differences in clustering coefficient and global efficiency distributions for individual networks (thin, black boxes) and the four classes of networks from clustering results (colored boxes). The majority of networks from each class have a clustering coefficient between 0 and 0.2 (Fig. \ref{fig:GS}a, top), with similar class median values. The networks show more variability in $E_{glob}$ (Fig. \ref{fig:GS}a, bottom), although similarly class median values are comparable. Local efficiency, characteristic path length, and modularity results are shown in Fig. \ref{fig:SGS}. These five computed graph statistical measures do not clearly distinguish between the classes recovered from clustering of homological features.
		% Then, need to end this paragraph with a summary statement about what we have learned that also foreshadows the next section.}
		
	\subsection{Class Members and Traits}
	
%\begin{landscape}
\begin{figure}[p]
\centering
    \includegraphics[width=6in]{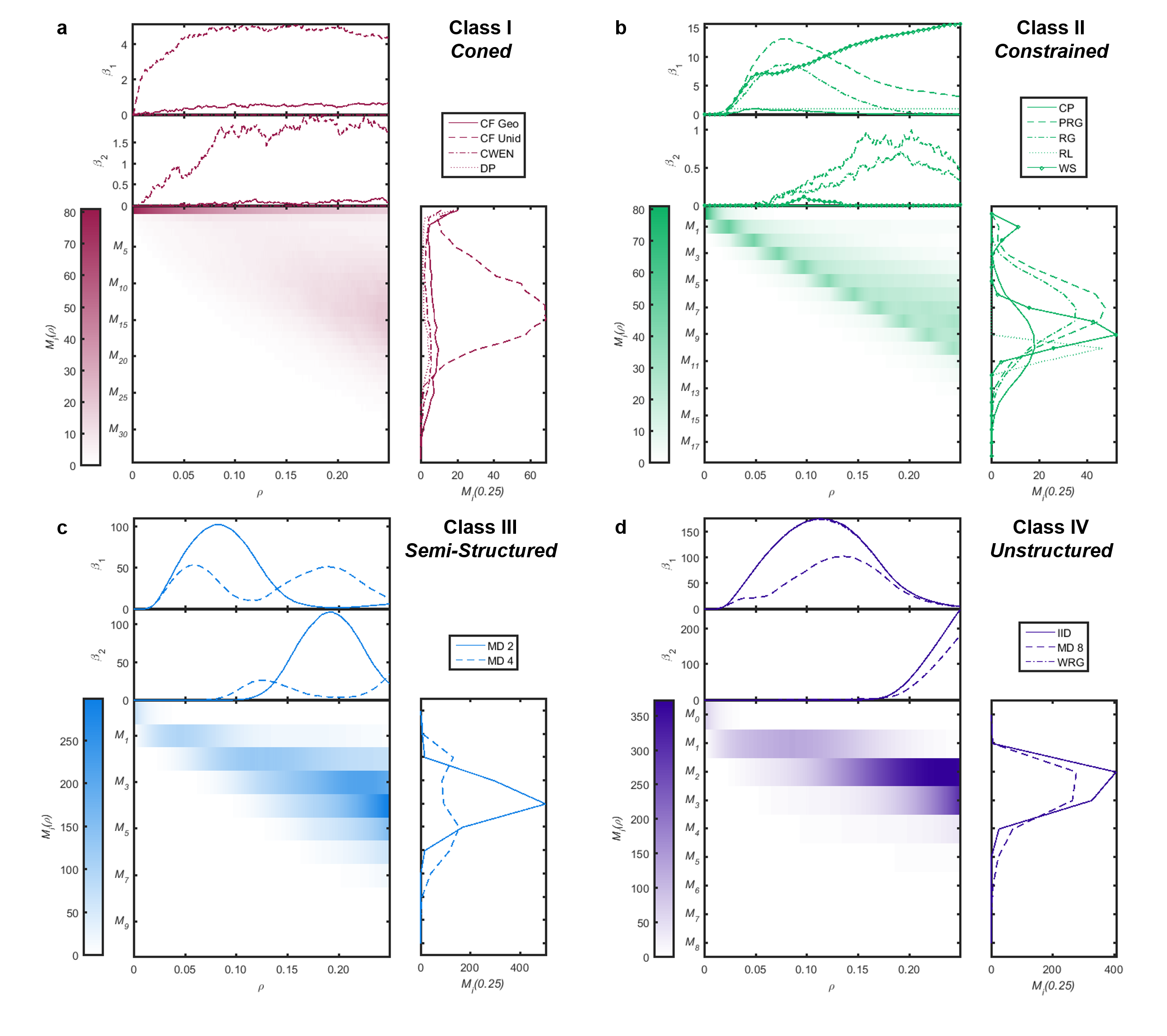}
    \caption{Model Network Clustering Reveals Four Topologically Defined Classes. Panels show the average network Betti curves for dimensions one and two ($\beta_1$, $\beta_2$; top), a class average heat map of the maximal clique evolution across edge density ($\rho$; left), and the network average maximal clique distribution at the final edge density ($\rho = 0.25$; right). \emph{(a)} Class I consists of configuration models from a geometric or discrete uniform distribution (CF Geo, CF Unid), the comprehensive weighted evolving network (CWEN), and the dot product (DP). \emph{(b)} Class II contains cross product (CP), probabilistic random geometric (PRG), random geometric (RG), ring lattice (RL), and Watts-Strogatz (WS) networks. \emph{(c)} Class III formed by modular networks with two or four communities (MD 2, MD 4). \emph{(d)} The independent and identically distributed (IID), weighted random graph (WRG), and modular network with eight communities (MD 8) form Class IV.}
    \label{fig:families}
\end{figure}
%\end{landscape}

	Next we want to understand the structural patterns that give rise to the four recovered classes. To examine such traits within each class, we compare Betti curves and maximal clique tracking results separated by class in Fig. \ref{fig:families} and integrate these with features from network construction. 
		
		Class I (Fig. \ref{fig:families}a) containing configuration models from a geometric or uniform distribution (CF Geo, CF Unid), the comprehensive weighted evolving network (CWEN) and the dot product (DP), shows the smallest amount of homology yet the largest maximal clique degree $\omega(G) = 35$. Generally $M(G)$ tracking is flat, due the high-dimensional cliques engulfing a large portion of the edges. The dot product (DP) graph is least similar, showing no homology and containing a large peak in the maximal clique distribution. 
		
		Members of this class contain a few high-strength, high-degree nodes that form cone-like structures with their neighbors, preventing cycle longevity or even cycle existence. This does not imply that the strength distributions of Class I graphs are scale-free or even similar, as clearly the CF Unid model has many more high-strength nodes than the other graphs. Instead, the cone structure speaks to the organization of links near particular cone points that may have high centrality\footnote{Recall centrality measures the number of shortest paths between node pairs that travel through the node of interest.}, but this is not required as cone points may not be a central point of the entire network. Unlike the other networks tested, Class I networks have a large $\overline{\beta}_{0}$ value. Recall $\overline{\beta}_{0}$ is the number of connected components, so we see a giant component emerge and smaller components that remain for many more edge additions than in other classes. We name Class I the \textit{Coned} group denoting the structural theme in these networks.
		
		Class II contains commonly studied networks such as the random geometric (RG), ring lattice (RL), and Watts-Strogatz (WS) models in addition to cross product (CP) and probabilistic random geometric models (PRG) (Fig. \ref{fig:families}b). These networks contain more homology than Class I networks, while also carrying a smaller $\omega(G)$. Networks in Class II are all generated from a sort of minimal-wiring technique. Network members strike a balance between high-dimensional cliques and number of cycles compared to Class I and IV which are skewed in one direction. This balance speaks to the geometric constraints imposed by minimal wiring, such as the triangle inequality. For any three nodes with two of the three pairs connected, the waiting period for the final edge to arrive is limited since its length, or the distance between end nodes, can be no longer than the sum of lengths between the connected nodes. In contrast, for a random graph the edge in question can continue to evade existence without any bound. Such constraints create networks where one might imagine their global structure as connecting many sizable clusters, occasionally cyclically. However, there is still a random component and these networks are not completely predefined, thus we call networks in this group \textit{Constrained}. 
		
		The third class contains only two networks, the modular network with two or four communities (MD 2, MD 4; Fig. \ref{fig:families}c). Both networks exhibit bimodal homology and asymmetric maximal clique distribution peaks. From these results it is reasonable to conclude that these form structures between that of Class II and IV. Furthermore, the construction of these modular networks is locally the same as with the weighted random graph, but globally they are two (four) small random networks tethered to each other by intermodule connections. Additionally, it is crucial that the ratio of community size to number of nodes is relatively large, as the modular network with eight communities is not included in this class. These features not only imposes a global structure but also allows smaller subnetworks to keep their random nature. Therefore, we call this the \textit{Semi-Structured} group.
		
		Finally, independent and identically distributed (IID), modular network with eight communities (MD 8), and the weighted random graph (WRG) comprise of Class IV (Fig. \ref{fig:families}d). These networks are very similar to each other in both Betti Curves and maximal clique distribution, particularly the IID and WRG. The IID and WRG networks are random by construction, with no constraints persuading them to form high dimensional cliques. The MD 8 graph shows a few higher dimensional cliques, but the small community size compared to the total number of nodes does not allow many internal high dimensional cliques to form, and therefore clique size, and subsequently global structure, is generally driven by the inter-module edges. 
		
		Each of the four recovered classes have a particular arrangement of architecture, seen from the number of cliques and cycles in the network. Furthermore, we see members within classes exhibit similar constructs of mesoscale homological features.

\subsection{Testing Example Networks}

With the recovered classes of weighted networks, we can now ask about the similarity in mesoscale structure between networks from biology and the presented models. We calculated the topological statistics of the networks formed from structural neural data (DSI), correlation between Kuramoto Oscillators with two or four communities (KM 2, KM 4), and protein-protein interactions in the citric acid cycle (TCA).

\begin{figure}[h!]
			\centering
				\includegraphics[width=6in]{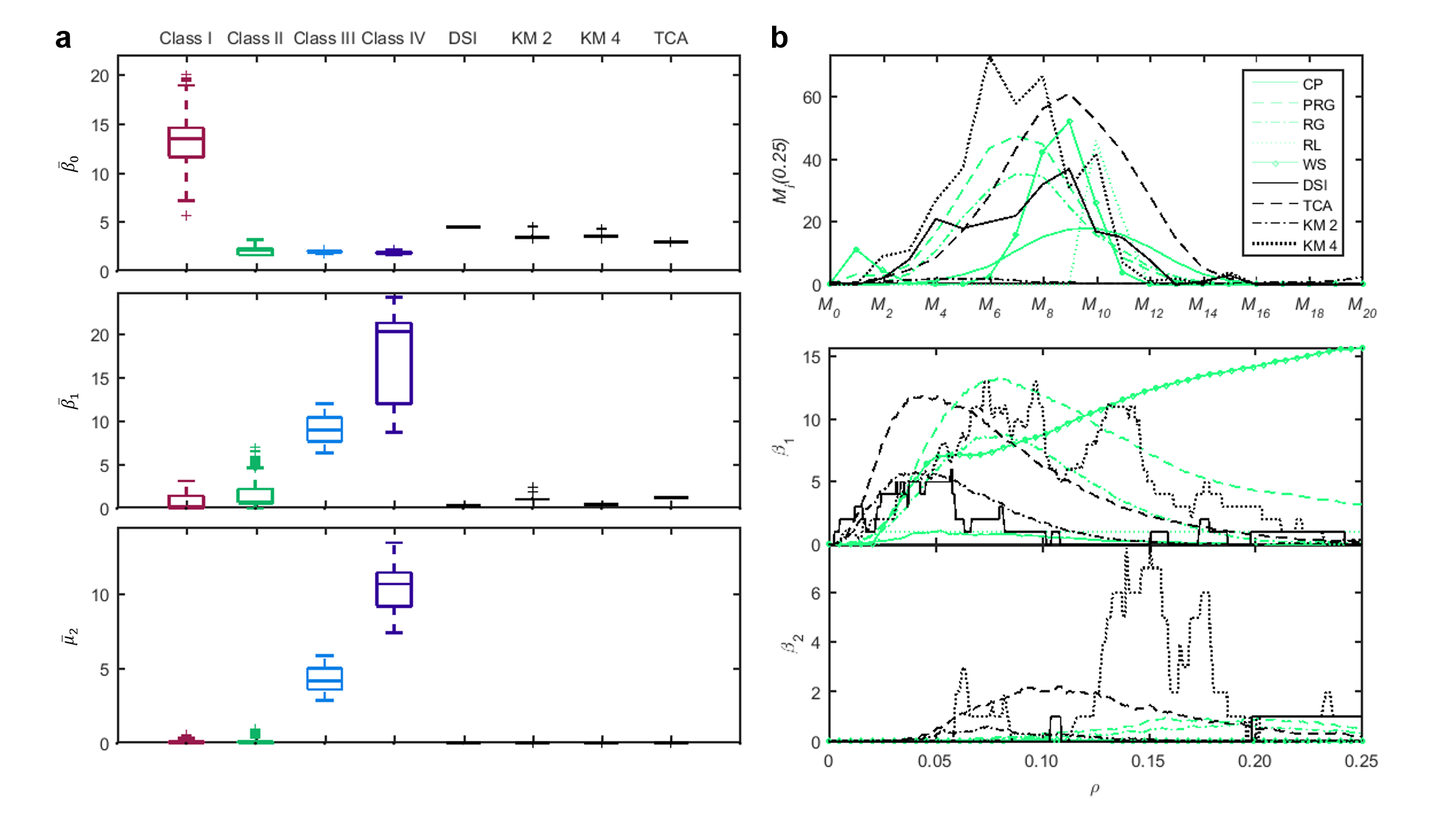}
				\caption{Biologically Inspired Networks Cluster with Constrained Class. \emph{(a)} Box plots showing distributions for recovered classes indicated by outline color and example networks in black, solid boxes. \emph{(b)} Maximal clique distribution at $\rho = 0.25$. DSI, TCA, and KM networks shown as black solid, dotted, and dashed lines, respectively. Class II members shown in pale green. Panels \emph{(b)} and \emph{(c)} share plot legends. \emph{(b)} Maximal clique distribution at $\rho = 0.25$. \emph{(c)} (Top) Betti curves $\beta_1$ and $\beta_2$ plotted across edge density ($\rho$).}
			\label{fig:clusterapp} 
		
	\end{figure}
Using distances from these values to class centroids (Fig. \ref{fig:centroid}a), these are most similar to the \emph{Constrained} class (Fig. \ref{fig:clusterapp}). By inspection, we see the distribution of cliques (Fig. \ref{fig:clusterapp}a) at $\rho = 0.25$ peak is similarly located, and the amount of homology seen in the Betti Curves (Fig. \ref{fig:clusterapp}c) and $\overline{\beta}$ and $\overline{\mu}$ values is comparable. Based on these results we conclude that these networks from biological systems and models are structurally most similar to networks in the \emph{Constrained} class.

%%%%%%%%%%%%%%%%%%%% DISCUSSION %%%%%%%%%%%%%%

%%%%%%%%%%%%%%%%%%%% DISCUSSION %%%%%%%%%%%%%%

\section{Discussion}

We have shown the homological features identified using persistent homology classify the 14 tested network models into four classes. This offers an alternative perspective on network structural similarity from that given by common graph statistics. Three tested biologically inspired networks display mesoscale homological features most similar to those from the cross product, random geometric, ring lattice, Watts-Strogatz, and probabilistic random geometric model networks.

\subsection{Classes Identified from Clustering}

%\ann{From DB: Describe what features homological statistics are sensitive to, and how this differs from other approaches. Then describe how these specific sensitivities explain the groupings into network model classes}

Homology detects mesoscale structure in the form of cliques and cycles in all dimensions. In our weighted networks, persistent homology sews together this information as the edge density parameter varies, assembling a sequence of blueprints from which we discern global architecture. This method is therefore sensitive to changes in both number and time of appearance of these features while remaining stable in the presence of minor edge reorderings \cite{cohen2007stability}. It is particularly important here to record both cliques and cycles, as cycles are ``mesoscale'' features while cliques are finer ``neighborhood'' features.  In contrast, many common graph statistics, such as those presented here, are aggregates of ``vertex-local'' measurements focused in low dimensions.

% ``hard'' and ``soft'' language from Persistent Homology and Shape Description by Carlsson et al. 2003. 

This feature duality is crucial here for capturing the similarity between Class I models as well as the differences in modular from random networks. Specifically, Class I models have varying strength distributions (from scale-free to sampled from a uniform distribution) yet carry a similar number of cycles. Additionally, the modular networks have local properties similar to the weighted random graph, but are distinct in both cycle time of birth and maximal clique count (``neighborhood'' and ``mesoscale'' features), thus distinguishing modular networks with a comparable community size to network ratio as their own class. Overall we see recovered classes contain similar topological values but within classes members can differ in topological signatures.

\subsection{Example Networks}

While prior studies of mesoscale community structure determined that protein-protein interaction and neural networks are distant cousins in the same family \cite{onnela2012taxonomies}, the homological mesoscale features computed here suggest that they are structurally quite similar. Indeed, the Betti curves and maximal clique distribution of the DSI, TCA, and KM models match those of the \emph{Constrained} class very closely (Fig. \ref{fig:clusterapp}). These three networks relate directly or indirectly to some functional task, and noting the maximal clique distribution peak we might envision these networks working with functional units of average size approximately six. We can see this as protein complex formation in the TCA network, or small groups of Kuramoto oscillators that highly correlate with one another. Having small modules of moderate size is perhaps most intriguing in the DSI data, however, as it may inform our understanding of brain function \cite{bassett2010efficient,bassett2011conserved}. Structural cliques formed by nearby brain regions may indicate possibilities for functional units \cite{chen2008revealing}. Moreover, these units may have dynamic functions that fluctuate based on brain task or state \cite{braun2015dynamic}, properties allowed by the mesoscale structural framework \cite{hermundstad2013,hermundstad2014,goni2014,honey2009}. 

The clique distribution of DSI and TCA networks is interestingly very similar, and both biological networks are closer in centroid distance to the CP, PRG, and RG networks than both RL and WS models (Fig. \ref{fig:centroid}b). While certainly the constraints imposed on these networks are not identical, both are under pressure to conserve energy \cite{bullmore2012,niven2008}. In the white matter tracts, this is more structural, as it is costly to lay excessive axonal projections \cite{bassett2010efficient}. In the protein-protein interaction network, it is more a conservation of energy in a metabolic sense \cite{zhang2009metabolic}. Indeed, though we should not immediately assume these interacting proteins are all directly along a metabolic pathway (some may be regulators, receptors, etc.), it is reasonable to consider the pressure on the cell to conserve energy in translating, transcribing, and maintaining health of proteins. Thus, this metabolic energy conservation may not only be realized in the pathways in which the queried proteins participate, but also in other processes that contribute to protein production and maintenance. In contrast to the brain network, the TCA network is a subset of a geometrically unbounded protein-protein interaction network of the entire proteome. One might question what we can infer about the entire network from work on this subset, which would be an interesting future avenue of research.

However, as seen clearly in the $\beta$ curves (Fig. \ref{fig:Sbettis2}), the construction processes of the DSI and TCA networks are distinct, implying that the backbones of these networks are wired differently.  Finally, the DSI network has a relatively high $\overline{\beta}_0$ value, indicating the presence of longer lived isolated components than those present in the other example networks \cite{stolz2014computational}. Indeed, the presence of these long-lived components is consistent with previous literature describing putative hub nodes in the large-scale brain networks \cite{van2011rich,hagmann2008mapping}.

\subsection{Comments on the Space of All Networks}

\begin{figure}[!h]
\begin{centering}
 \includegraphics[width=6in]{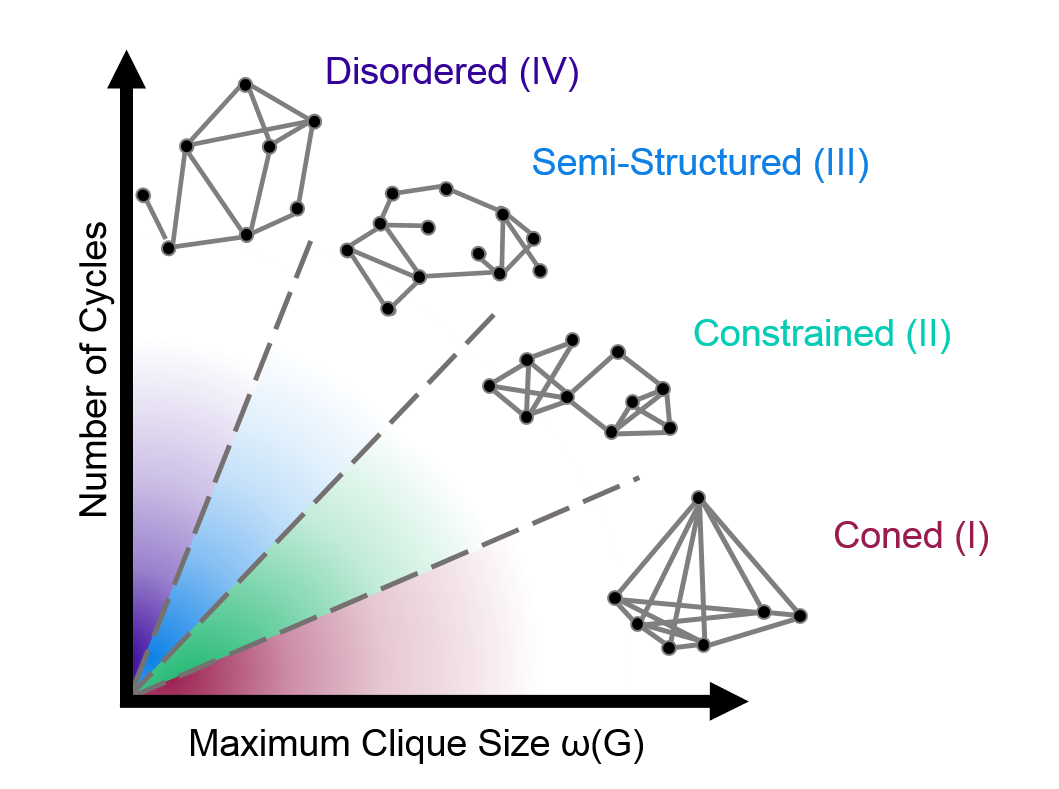}
\caption{Clustering Results Offer Conceptual Organization on the Space of Possible Networks.  Considering all possible networks on a given number of nodes, the four classes found (Coned, Constrained, Semi-Structured, and Disordered) span regions of this schematic plane encompassing networks with similar ratios of number of cycles to maximum clique size $\omega(G)$.}
\label{fig:disc} 
\end{centering}
\end{figure}

While we use clusters as a conceptual tool to understand common network models, it is fair to consider a spectrum of potential network structures (Fig. \ref{fig:disc}). At one end, the \emph{Unstructured} group has many low degree cliques and contains a larger number of cycles. Indeed, IID networks have been well characterized \cite{kahle2011random} and we expect a large number of small cycles to be present. As $\omega(G)$ increases, we often see fewer cycles and more higher-dimensional cliques, a finding that is similar to observations in the \emph{Semi-Structured} and then \emph{Constrained} class. In networks with high $\omega(G)$ we often find a few nodes with a high relative strength forming cone points, killing much of the homology. Extrapolating to the extreme case gives a network that is one maximal clique.

\subsection{Further Applications of Persistent Homology}

	In this paper, we have used some summary outputs of the persistent homology computation. However, more statistics can in principle be extracted, and might offer additional insights into network structure. For example, recent advances make it possible to recover minimal generators for the homology \cite{henselmannovel}. Particularly useful for empirical networks, identifying generators can show the nodes responsible for these interesting structural features. In ongoing research, we use these tools to uncover mesoscale architecture in the structural neural network. Seen here in the TCA network, the nonzero $\beta_3$ curve indicates the presence of 3-cycles. Knowing the enzymes responsible for this structure may inform metabolic network efficiency and robustness. Indeed, particular mutations may create or disrupt cycles, offering information about the structural role of the wild-type protein.  We speculate that cycle generators in neural networks could inform brain communication and function \cite{stolz2014computational}. Cycle participants or even cycle presence may differ based on age, condition, or previous training. Presence alone may indicate a pathway that diverges and reconverges, such as the dorsal/ventral visual streams \cite{sepulcre2012stepwise}, or perhaps the possibility of information flow circumventing hub nodes, which will be the aim of future studies.

\section{Conclusions}
Persistent homology is sensitive to the global structure of a weighted network, making it an effective tool for network classification. The power of these topological methods is clearly demonstrated in the identification of four network classes that are indistinguishable using vertex-centric graph statistics, emphasizing the importance of considering both the global and local structure of a weighted network. Moreover, we observe that networks from DSI data, protein-protein interactions, and coupled Kuramoto oscillators are structurally most similar to minimally wired networks. These data argue for the importance of considering a topological perspective in understanding network structure in real and synthetic networked systems.

\section{Acknowledgments}
DSI data was collected and preprocessed by Matthew Cieslak and Scott Grafton, in the Department of Psychological and Brain Sciences at the University of California, Santa Barbara. DSB acknowledges support from the John D. and Catherine T. MacArthur Foundation, the Alfred P. Sloan Foundation, the Army Research Laboratory and the Army Research Office through contract numbers W911NF-10-2-0022 and W911NF-14-1-0679,the National Institute of Mental Health (2-R01-DC-009209-11), the National Institute of Child Health and Human Development (1R01HD086888-01), the Office of Naval Research, and the National Science Foundation (BCS-1441502 and BCS-1430087). CG acknowledges support from the Warren Center for Network and Data Sciences.


\bibliographystyle{plain}
\bibliography{bibfile}

\newpage
\appendix
\section{Appendix}

%\DBadvice{Throughout the appendix, make sure that the equations are appropriately punctuated, either being followed by a period or a comma. Also add original citations for the first definitions of these statistics. If you need help locating them, let me know -- I suggest looking in Phys Rev E Stat Nonlin Soft Matter Phys. 2012 Oct;86(4 Pt 1):041306.}

%\counterwithin{figure}{section}
\renewcommand{\thefigure}{A\arabic{figure}}

\setcounter{figure}{0}

\subsection{Clustering Results}

\begin{figure}[h!]
\begin{centering}
 \includegraphics[width=4in]{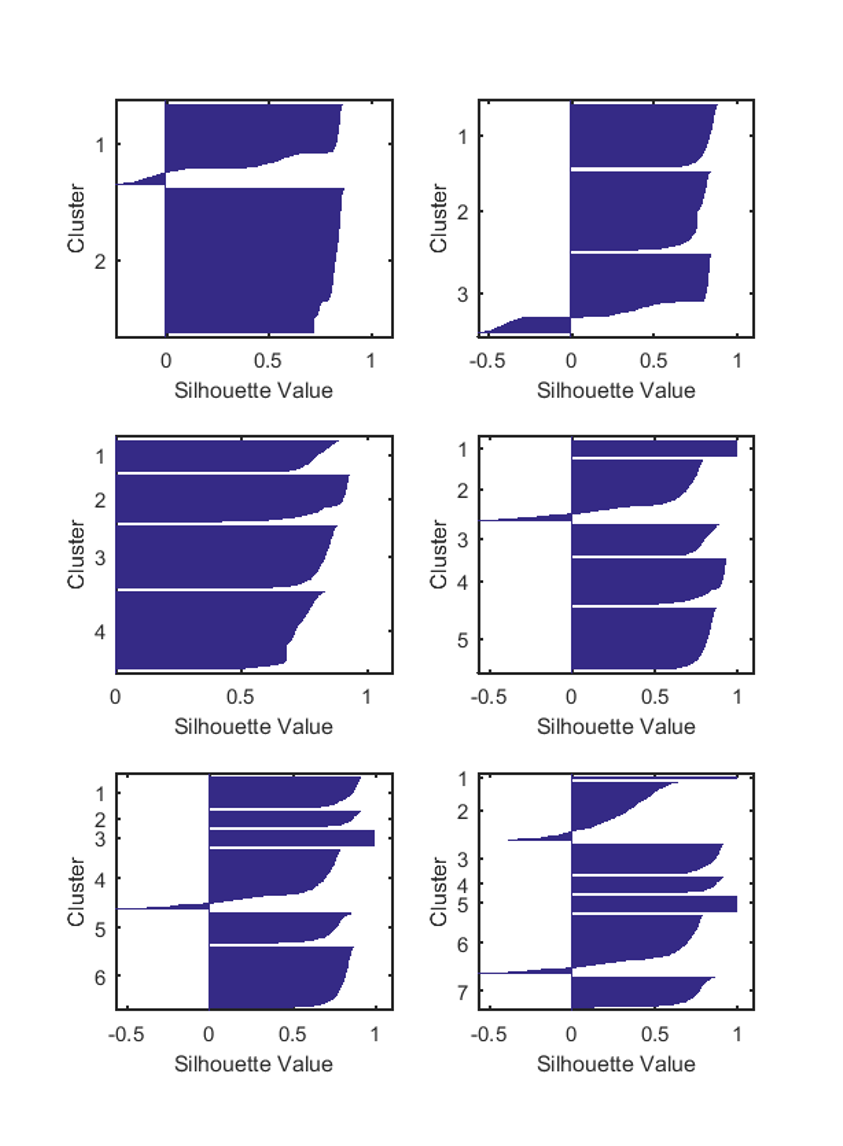}
\caption{Silhouette Plots from Hierarchical Clustering Suggest Four Classes of Model Networks. Plots show the similarity of networks to members of the same class in comparison to outside networks (horizontal axis). The number of resulting clusters varies from two (top left) to seven (bottom right).}
\label{fig:silplot} 
\end{centering}
\end{figure}

To determine the optimal number of classes, we analyzed silhouette plots on clustering results with different number of classes specified (Fig. \ref{fig:silplot}).  

\begin{figure}[h!]
\begin{centering}
 \includegraphics[width=5in]{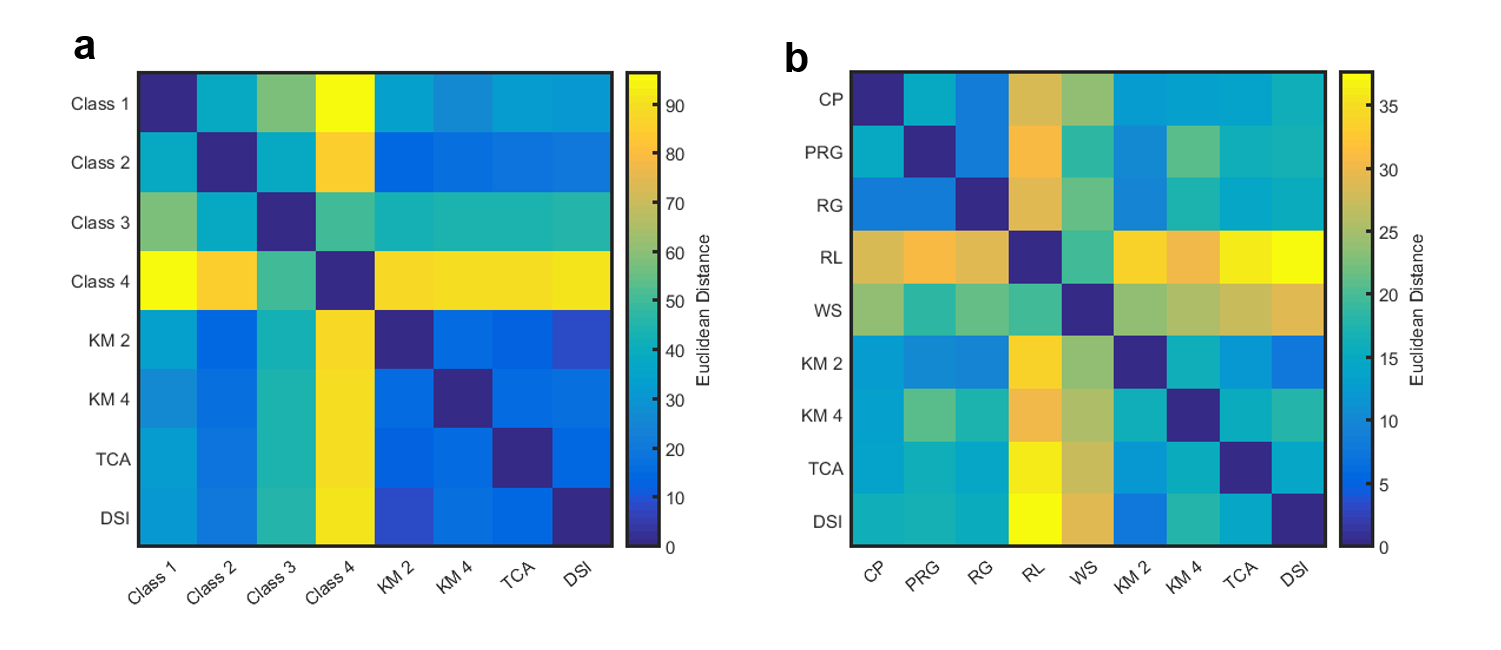}
\caption{Similarity of Example Networks to Constrained Class Models. (\emph{a}) Distance matrix indicating Euclidean distance between centroids of recovered classes and biologically inspired networks. (\emph{b}) Euclidean distances between centroids of Class II models and example networks.}
\label{fig:centroid} 
\end{centering}
\end{figure}

For simple determination of the most similar class to biologically inspired networks, we computed the Euclidean distance between the four class centroids and the four example networks (Fig. \ref{fig:centroid}a). Within the Constrained class, Fig. \ref{fig:centroid}b shows the distance between centroids of individual models and biologically inspired networks (Fig. \ref{fig:centroid}b). 

\subsection{Graph Statistics}
	\textit{Clustering Coefficient ($C$)} For each node, we calculated the weighted clustering coefficient 
	
\begin{align}\label{CC}
\begin{split}
 C_i = \frac{1}{k_i(k_i - 1)} \sum_{j,k \in V}(\hat{w}_{i,j}\hat{w}_{i,k}\hat{w}_{j,k})^{1/3}
\end{split}
\end{align}
 where $\hat{w}_{i,j} = w_{i,j}/$max$(w_{i,j})$ for $i,j,k \in V$ \cite{rubinov2010complex}. The overall clustering coefficient of the graph ($C$) is then the average of the individual clustering coefficients of all nodes.
	 
	 \textit{Global Efficiency ($E_{glob}$)} Using the inverse of the weighted distance $d_{i,j}^W$ between nodes $i$ and $j$ we calculated the global efficiency, defined as 
	 	
	 	\begin{align}\label{Eglob}
\begin{split}
	 	E_{glob,i} = \frac{1}{n-1}\sum_{j \in N, i \neq j}(d_{ij}^W)^{-1}.
	 	\end{split}
\end{align}
The overall global efficiency $E_{glob}$ of the network is the mean of $E_{glob,i}$ over all nodes \cite{rubinov2010complex}. 
	 
	 \textit{Local Efficiency ($E_{loc}$)} At every node $i$ the local efficiency is defined 
	 
	 \begin{align}\label{Eloc}
\begin{split}
	 E_{loc,i} = \frac{\sum_{j,h \in N, j\neq i}(w_{i,j}w_{i,h}[d_{jh}^W(N_i)]^{-1})^{1/3}}{k_i(k_i-1)}
	 \end{split}
\end{align}
 where $d_{jh}^W(N_i)$ is the inverse shortest path length for neighbors $j$ and $h$ of node $i$ in the set $N_i$ which includes all nodes except $i$. Taking the average of $E_{loc,i}$ over all $i \in N$ gives the local efficiency of the network, $E_{loc}$ \cite{rubinov2010complex}.
	 
	 \textit{Characteristic Path Length ($L$)} The characteristic path length is an average of all shortest paths between nodes within a connected component and is defined as
	 
	 \begin{align}\label{CPL}
\begin{split}
	 L_i = \frac{\sum_{j \in N, j \neq i}d_{ij}^W}{n-1}.
	 \end{split}
\end{align}
The characteristic path length $L$ is the average of $L_i$ over all nodes \cite{rubinov2010complex}.\\
	 
	 \textit{Modularity ($Q$)} Using the Louvain algorithm for community detection, we calculate the modularity of a network
	 
	 \begin{align}\label{Q}
\begin{split}
	   Q = \frac{1}{v} \sum_{i,j}(w_{i,j} - \frac{s_i s_j}{v})\delta_{M_i M_j}
	   \end{split}
\end{align}
 with the sum of all connection weights $v = \sum_{i,j}w_{i,j}$,  the strength of node $i$ defined as $s_i = \sum_j w_{i,j}$, and $M_i$ the community of node $i$. The modularity is the average difference between the weight of the intermodule connection and the expected weight \cite{rubinov2011weight}. 
 
 \begin{figure}[h!]
\begin{centering}
 \includegraphics[width=4in]{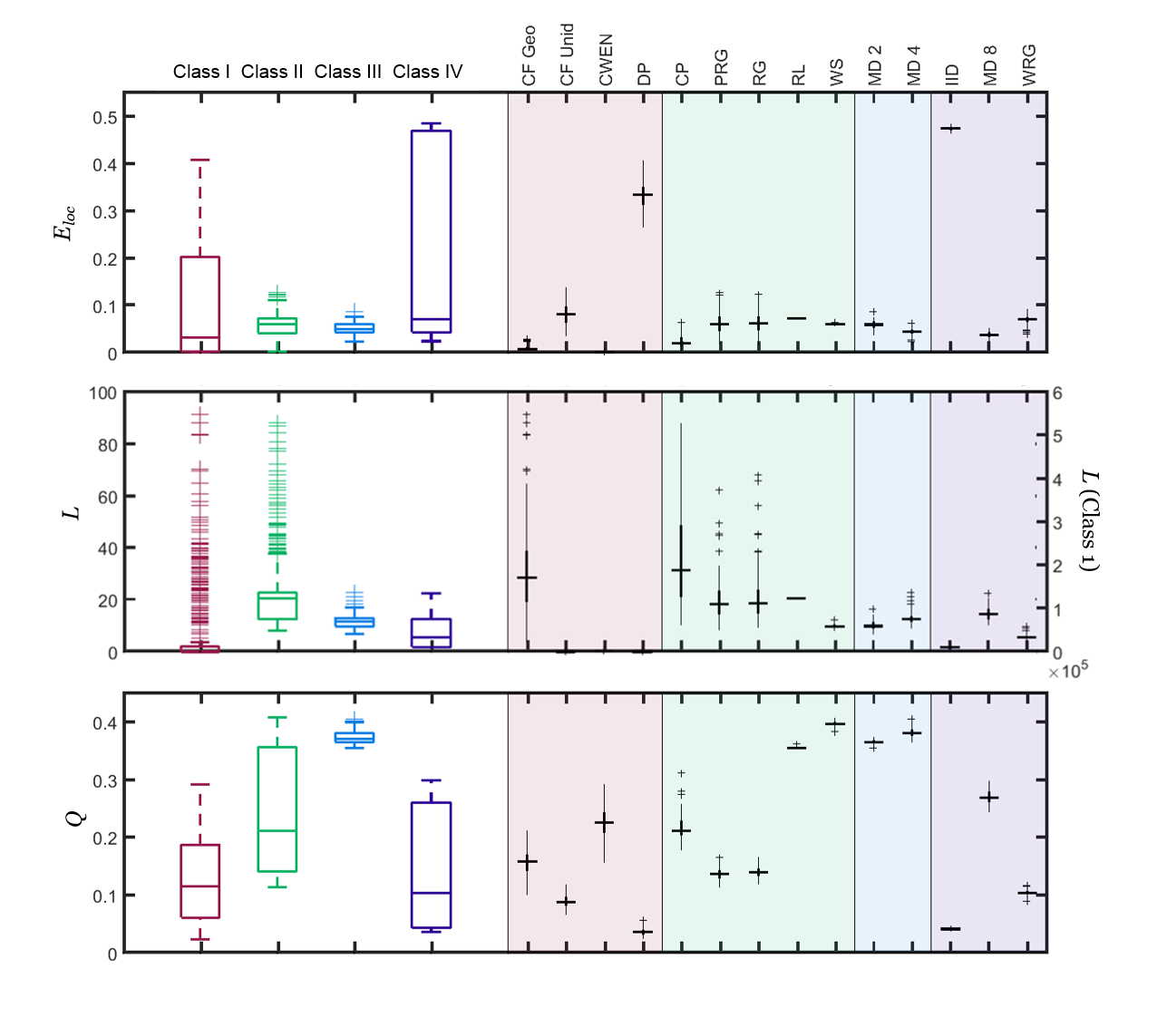}
\caption{Multiple Classes Share Similar Graph Metrics. Local efficiency ($E_{loc}$) (top), modularity ($Q$) (middle), and characteristic path length ($L$)  (bottom) of model networks, sorted by recovered groups. Each statistic shown by group (left, outlined with class color) and individually (right, black filled boxes shaded by color). Vertical axis of characteristic path length shown on the right for Class 1 networks individually and as a class.} 
\label{fig:SGS}
\end{centering}
\end{figure}

\subsection{Homological Features}

\begin{figure}[!h]
\begin{centering}
 \includegraphics[width=4.3in]{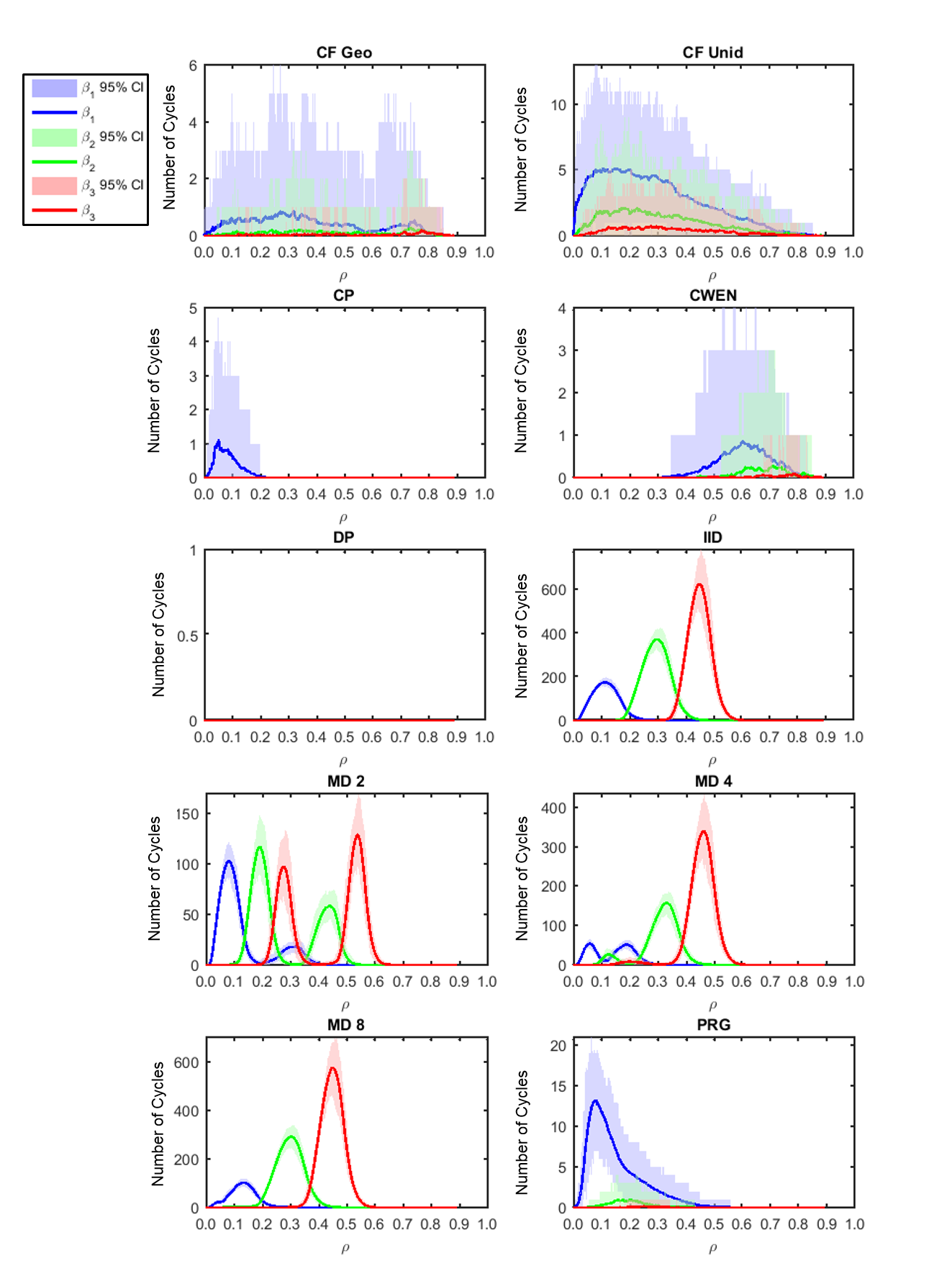}
\caption{Betti Curves Record Cycles Throughout the Weighted Network Filtration. Betti curves ($\beta_d$) in dimensions one (blue), two (green), and three (red) with 95\% confidence intervals indicated by shaded region. Configuration models from a geometric or discrete uniform strength distribution (CF Unid, CF Geo), cross product (CP), Comprehensive Weighted Evolving Network (CWEN), dot product (DP), independent and identically distributed (IID), modular with two (MD 2), four (MD 4), and eight (MD 8) communities, and probabilistic random geometric (PRG) networks shown.}
\label{fig:Sbettis1} 
\end{centering}
\end{figure}

\begin{figure}[!h]
\begin{centering}
 \includegraphics[width=4.3in]{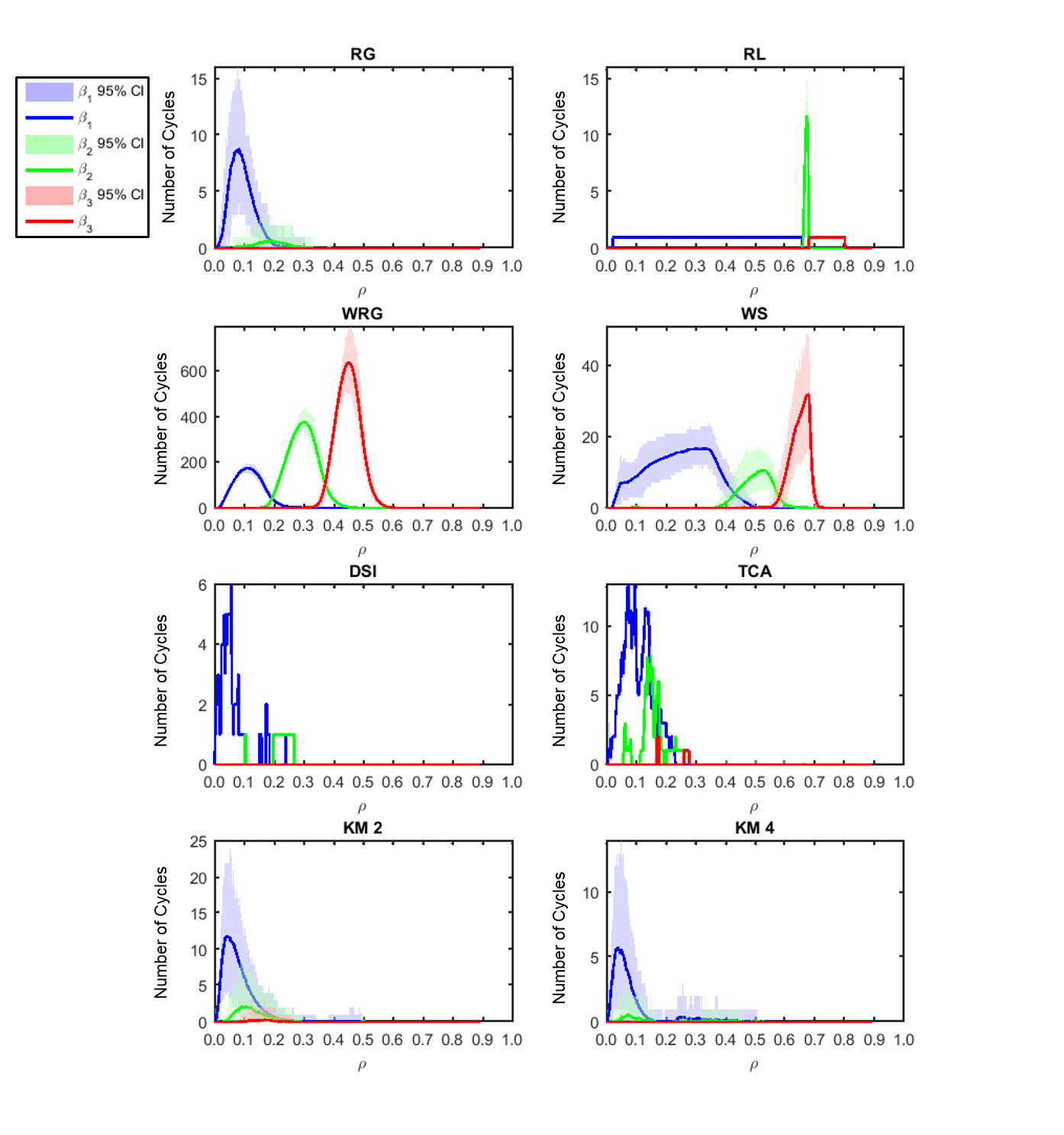}
\caption{Betti Curves Record Cycles Throughout the Weighted Network Filtration (Continued). Betti Curves ($\beta_d$) in dimensions one (blue), two (green), and three (red) with 95\% confidence intervals indicated by shaded region. Random geometric (RG), ring lattice (RL), weighted random graph (WRG), Watts-Strogatz (WS), neural structural connectome (DSI), protein-protein interaction in the tricarboxylic acid cycle (TCA), and Kuramoto Oscillator coupling networks with two (KM 2) and four (KM 4) communities shown.}
\label{fig:Sbettis2} 
\end{centering}
\end{figure}

To understand progression of cycles as edge weight increases, we recorded the Betti curves for each network model (Fig. \ref{fig:Sbettis1}, \ref{fig:Sbettis2}).

\begin{figure}[!h]
\begin{centering}
 \includegraphics[width=5in]{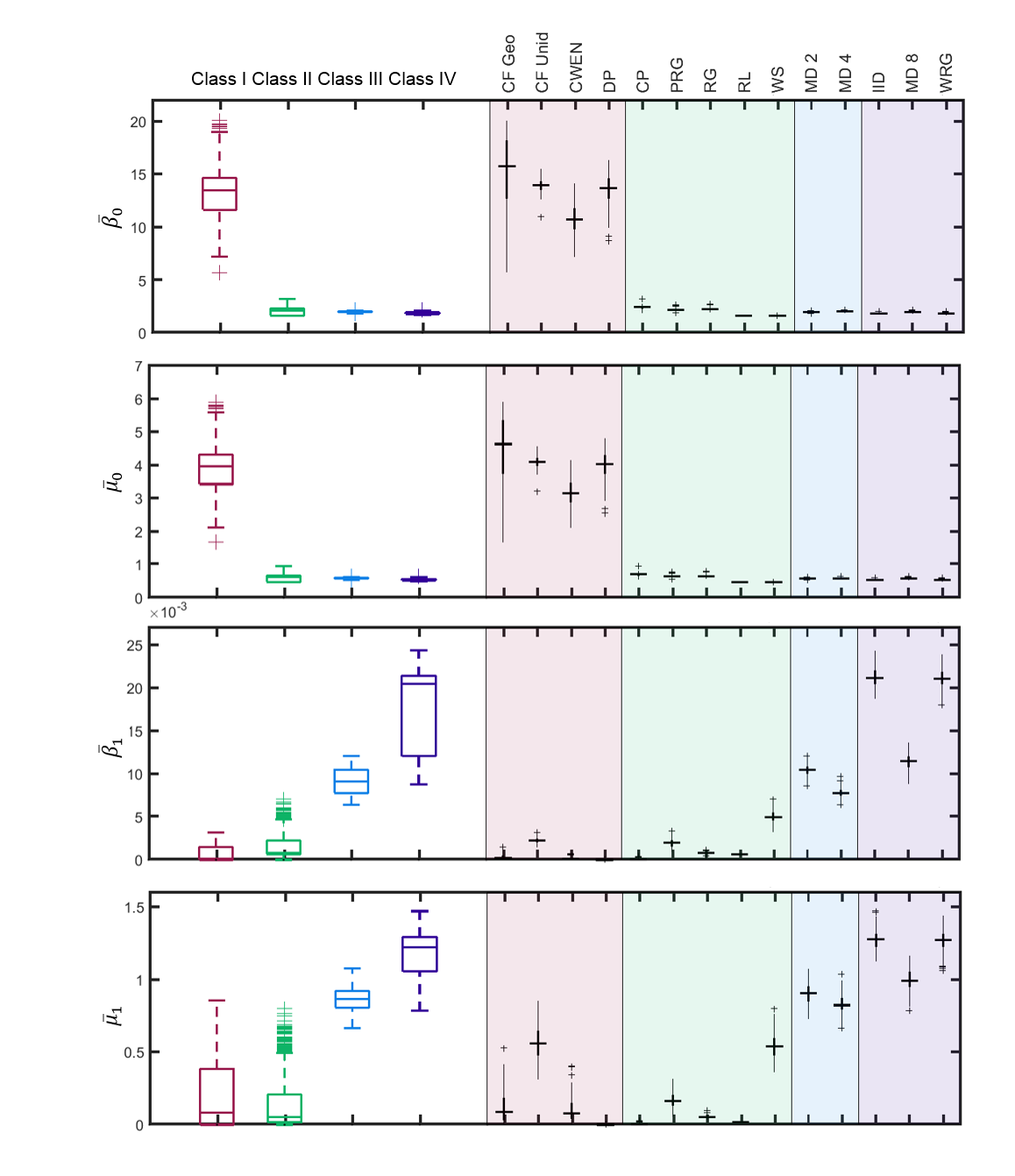}
\caption{Topological Statistics in Dimensions Zero and One Differ Between Recovered Network Classes. Calculated $\overline{\beta}_0$ (top), $\overline{\mu}_0$, $\overline{\beta}_1$. and $\overline{\mu}_1$ (bottom) values for recovered class (left, box outlined with class color) and model networks (right, shaded by class color).}
\label{fig:Sbbar1} 
\end{centering}
\end{figure}

\begin{figure}[!h]
\begin{centering}
 \includegraphics[width=5in]{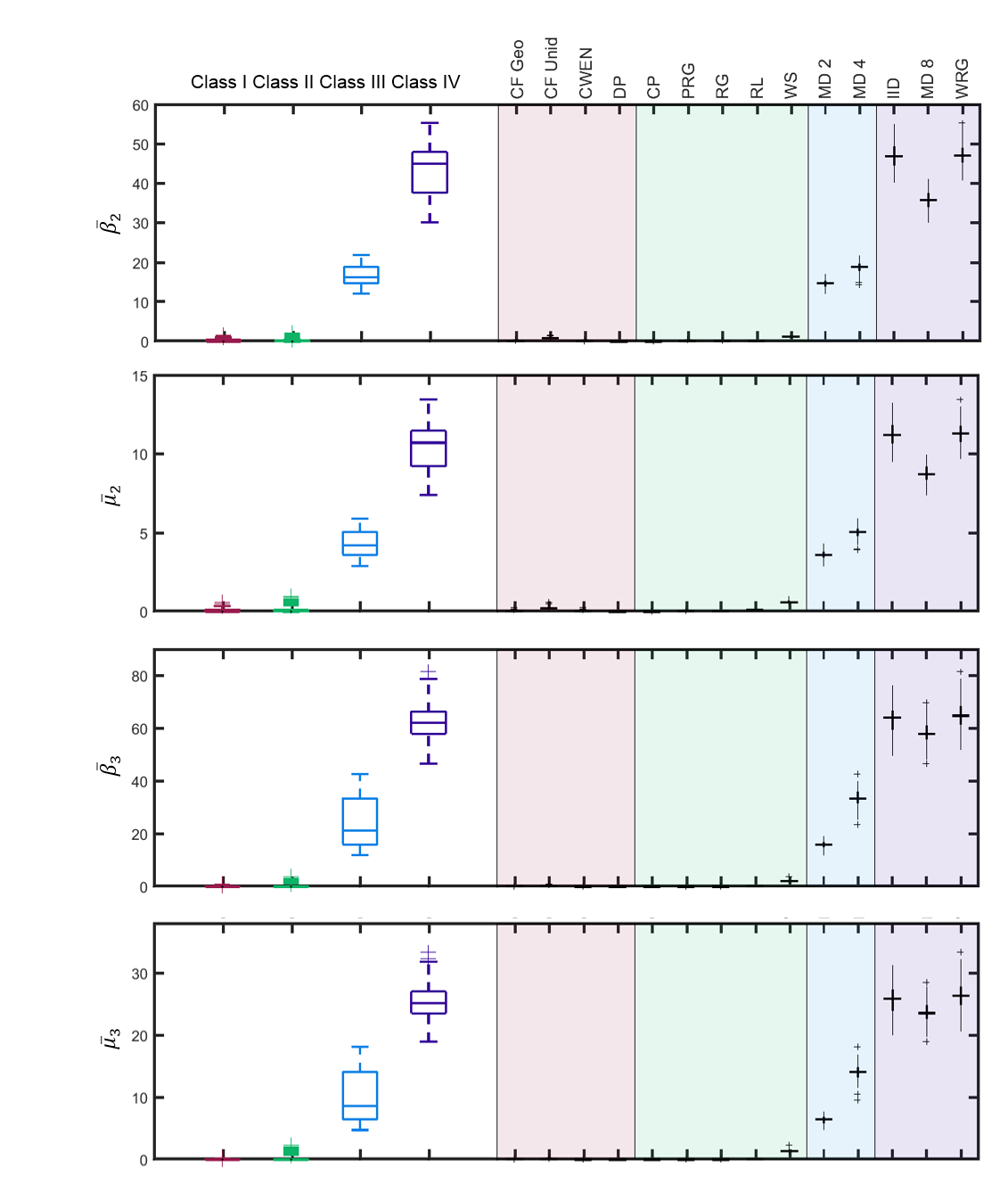}
\caption{Topological Statistics in Dimensions Two and Three Further Distinguish Classes. Calculated $\overline{\beta}_2
$, $\overline{\mu}_2$, $\overline{\beta}_3$. and $\overline{\mu}_3$ (bottom) values for recovered class (left, box outlined with class color) and model networks (right, shaded by class color).}
\label{fig:Sbbar2} 
\end{centering}
\end{figure}

\begin{figure}[!h]
\begin{centering}
 \includegraphics[width=5in]{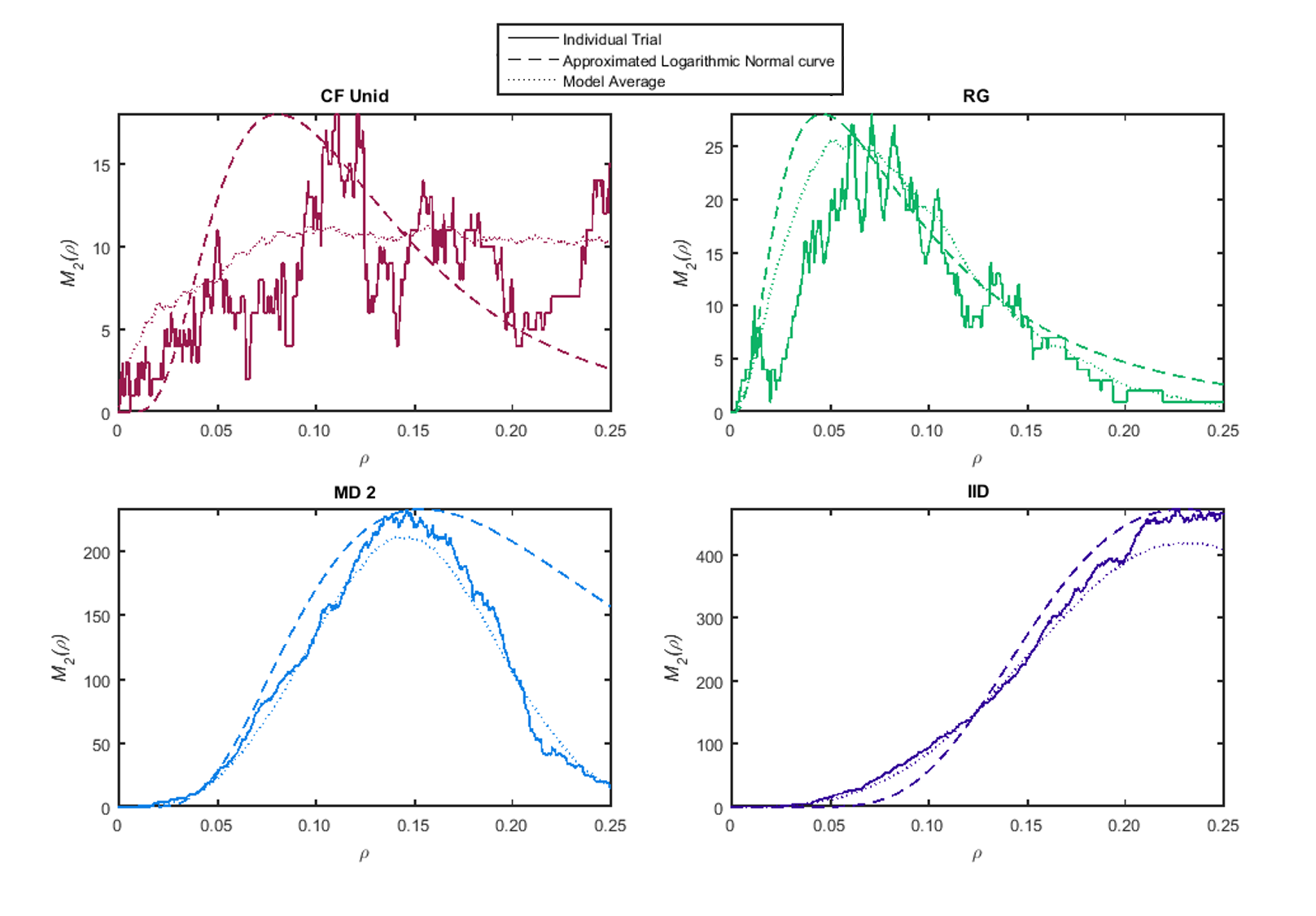}
\caption{Logarithmic Normal Distribution Approximates Maximal Clique Vector. In each plot, $M_3(\rho)$ for an individual model run (solid line), the scaled approximated logarithmic normal distribution with parameters from this run (dashed line), and the average of all runs within indicated model (dotted line) are shown across edge density ($\rho$). One model from each class is presented: configuration model with discrete uniform distribution (CF Unid; top left), random geometric (RG; top right), modular with two communities (MD 2; bottom left), and independent and identically distributed (IID; bottom right).}
\label{fig:sfit} 
\end{centering}
\end{figure}

For all runs in each model, we calculated the $\overline{\beta}_d$ and $\overline{\mu}_d$ values for $d = 0,1,2,3$ (Fig. \ref{fig:Sbbar1},\ref{fig:Sbbar2}) as well as approximated a logarithmic normal curve to the maximal clique vector $M$ using $\mu = \text{mean}(\ln(M_k))$ and $\sigma = \text{std}(\ln(M_k))$ (Fig. \ref{fig:sfit}) \cite{mood1974introduction}. Hierarchical clustering as described in Methods was performed with these features.

%\bibliography{sample}
%
% once the .bbl file has been generated then place the text in your article.

% To get the unnumbered reference style the author should use [unnumbib]
%as an option in the document class.  For example: \documentclass[unnumbib]{comnet}

\end{document}